\documentclass{article}%
\usepackage{graphicx}
\usepackage{amsmath}
\usepackage{amsfonts}
\usepackage{amssymb}%
\newtheorem{theorem}{Theorem}

\newtheorem{definition}[theorem]{Definition}

\begin{document}

\title{Geometric structures on\\$G_{2}$ and $Spin\left(  7\right)  $-manifolds}
\author{Jae-Hyouk Lee and Naichung Conan Leung}
\maketitle

\begin{abstract}
This article studies the geometry of moduli spaces of $G_{2}$-manifolds,
associative cycles, coassociative cycles and deformed Donaldson-Thomas
bundles. We introduce natural symmetric cubic tensors and differential forms
on these moduli spaces. They correspond to Yukawa couplings and correlation
functions in M-theory.

We expect that the Yukawa coupling characterizes (co-)associative fibrations
on these manifolds. We discuss the Fourier transformation along such
fibrations and the analog of the Strominger-Yau-Zaslow mirror conjecture for
$G_{2}$-manifolds.

We also discuss similar structures and transformations for $Spin\left(
7\right)  $-manifolds.  

\end{abstract}

\newpage

The \textit{mirror symmetry conjecture} for Calabi-Yau threefolds, which is
originated from the study of the string theory on $X\times\mathbb{R}^{3,1}$,
has attracted much attention in both the physics and mathematics communities.
Certain parts of the conjecture should hold true for Calabi-Yau manifolds of
any dimensions. But others are designed only for threefolds, for example, the
special K\"{a}hler geometry on the moduli space, the holomorphic Chern-Simons
theory and relationships to knot invariants. These features should be
interpreted in the realm of the geometry of the $G_{2}$-manifold, $M=X\times
S^{1}$. From a physical point of view, they come from the studies of M-theory
on $M\times\mathbb{R}^{3,1}$.

In mathematics, manifolds with $G_{2}$ holonomy group have been studied for a
while. In the Berger's list of holonomy groups of Riemannian manifolds, there
are a couple of exceptional holonomy groups, namely $G_{2}$ and $Spin\left(
7\right)  $. Manifolds with these holonomy groups are Einstein manifolds of
dimension seven and eight respectively. As cousins of K\"{a}hler manifolds,
the analogs of complex structures, complex submanifolds and
Hermitian-Yang-Mills bundles are vector cross products (Gray \cite{Gr}),
calibrated submanifolds (Harvey and Lawson \cite{HL}) and Donaldson-Thomas
bundles (\cite{DT}). These geometries are studied by many people in the Oxford
school including Donaldson, Hitchin, Joyce, Thomas and others.

A manifold $M$ with $G_{2}$-holonomy, or simply a $G_{2}$-manifold, can be
characterized by the existence of a parallel positive three form $\Omega$.
Because of the natural inclusion of Lie groups $SU\left(  3\right)  \subset
G_{2}$, the product of a Calabi-Yau threefold $X$ with a circle $S^{1}$ has a
canonical $G_{2}$-structure with $\Omega=\operatorname{Re}\Omega_{X}%
-\omega_{X}\wedge dt$, where $\Omega_{X}$ and $\omega_{X}$ are the holomorphic
volume form and the Calabi-Yau K\"{a}hler form on $X$.

The mirror symmetry conjecture for Calabi-Yau manifolds roughly says that
there is a duality transformation from the symplectic geometry (or the
A-model) to the complex geometry (or the B-model) between a pair of Calabi-Yau
threefolds. For instance natural cubic structures, called the A- and B-Yukawa
couplings, on their moduli spaces can be identified. This gives highly
non-trivial predictions for the enumerative geometry of $X$.

From physical considerations, it is more natural to understand the
\textit{M-theory} on $M\times\mathbb{R}^{3,1}$ with $M$ being a $G_{2}%
$-manifold. This is studied recently by Acharya, Atiyah, Vafa, Witten, Yau,
Zaslow and others (e.g. \cite{Ac}, \cite{AV}, \cite{AMV}, \cite{AW}). To
better understand the M-theory and its duality transformations, or
(equivalently) the geometry of $G_{2}$-manifolds, we need to study natural
geometric structures on various moduli spaces attached to $M$. For example we
introduce the analog of the Yukawa coupling on the moduli space of $G_{2}%
$-metrics as follows,
\[
\mathcal{Y}_{M}\left(  \phi\right)  =\int_{M}\Omega\left(  \hat{\phi}%
,\hat{\phi},\hat{\phi}\right)  \wedge\ast\Omega,
\]
where $\phi\in H^{3}\left(  M,\mathbb{R}\right)  $. This cubic structure comes
from the exceptional Jordan algebra structure for $E_{6}$ and it seems to be
even more natural than those A- and B-Yukawa couplings for Calabi-Yau
threefolds. It will be used later to characterize different kinds of fibration
structures on $M$.

\bigskip\ 

The geometry of a $G_{2}$-manifold is reflected by its calibrated submanifolds
(\cite{HL}) and Yang-Mills bundles (\cite{DT}). Calibrated submanifolds are
always volume minimizing and there are two different types in a $G_{2}%
$-manifold $M$: (i) coassociative submanifolds $C$ of dimension four,
calibrated by $\ast\Omega$ and (ii) associative submanifolds $A$ of dimension
three, calibrated by $\Omega$. Sometimes it is important to include both
points and the whole manifold in the above list, together they are calibrated
by different components of $e^{\Omega+\ast\Omega}$.

In physics, one needs to study \textit{supersymmetric cycles }\cite{MMMS}
which are calibrated submanifolds together with \textit{deformed Yang-Mills
bundles }over them. In the $G_{2}$ case, there are three different types:
\textit{Coassociative cycle }$\left(  C,D_{E}\right)  $, with $D_{E}$ an ASD
connection on a coassociative submanifold $C$; (ii) \textit{Associative cycle
}$\left(  A,D_{E}\right)  $, with $D_{E}$ a unitary flat connection on an
associative submanifold $A$, it is a critical point of the following
functional,
\[
\int_{A\times\left[  0,1\right]  }Tr\left[  e^{\ast\Omega+\tilde{F}}\right]
.
\]
(iii) \textit{Deformed Donaldson-Thomas connection }$D_{E}$ which satisfies
the equation
\[
F\wedge\ast\Omega+F^{3}/6=0,
\]
and it is a critical point of the following functional,
\[
\int_{M\times\left[  0,1\right]  }Tr\left[  e^{\ast\Omega+\tilde{F}}\right]
.
\]

Moduli spaces of these cycles are denoted as $\mathcal{M}^{coa}\left(
M\right)  ,\mathcal{M}^{ass}\left(  M\right)  $ and $\mathcal{M}^{bdl}\left(
M\right)  $ respectively, or simply $\mathcal{M}^{coa},\mathcal{M}^{ass}$ and
$\mathcal{M}^{bdl}$. We will define canonical three forms and four forms on
them using the Clifford multiplication on spinor bundles and the Lie algebra
structure on the space of self-dual two forms. In physics languages, they
correspond to \textit{correlation functions }in quantum field theory. We will
also introduce natural symmetric cubic tensors on these moduli spaces.

In the flat situation such canonical three forms determine $G_{2}$-structures
on both $\mathcal{M}^{coa}$ and $\mathcal{M}^{ass}$. In general the moduli
space of associative (resp. coassociative) submanifolds can be regarded as a
coassociative (resp. associative) subspace of $\mathcal{M}^{ass}$ (resp.
$\mathcal{M}^{coa}$).

We summarize these structures on various moduli spaces in following tables:

\begin{center}
$%
\begin{tabular}
[c]{|l|l|l|l|l|}\hline
& $\text{moduli }$ & Yukawa coupling$\text{ }$ & Metric tensor &
Prepotential\\\hline
$\mathcal{M}^{G_{2}}$ & $G_{2}\text{-structures}
\begin{array}
[c]{l}%
\,\\
\,
\end{array}
$ & \multicolumn{1}{|c|}{$\mathcal{Y}$} & \multicolumn{1}{|c|}{$\mathcal{G}$}
& \multicolumn{1}{|c|}{$\mathcal{F}$}\\\hline
\end{tabular}
$

{\small Table: Structures on the moduli space of G}$_{2}${\small -metric with
unit volume.}
\end{center}

\bigskip

\begin{center}%
\begin{tabular}
[c]{|l|l|l|l|l|}\hline
& $\text{moduli of}$ & $\text{3 form }$ & $\text{4 form}$ & $\text{cubic
tensor}$\\\hline
$M$ & point$\text{s}
\begin{array}
[c]{l}%
\,\\
\,
\end{array}
$ & \multicolumn{1}{|c|}{$\Omega$} & \multicolumn{1}{|c|}{$\Theta$} &
\multicolumn{1}{|c|}{}\\\hline
$\mathcal{M}^{coa}\left(  M\right)  $ & $%
\begin{array}
[c]{l}%
\text{coassociative }\\
\text{cycles}%
\end{array}
$ & \multicolumn{1}{|c|}{$\Omega_{\mathcal{M}^{coa}\left(  M\right)  }$} &
\multicolumn{1}{|c|}{$\Theta_{\mathcal{M}^{coa}\left(  M\right)  }$} &
\multicolumn{1}{|c|}{$\mathcal{Y}_{\mathcal{M}^{coa}\left(  M\right)  }$%
}\\\hline
$\mathcal{M}^{ass}\left(  M\right)  $ & $%
\begin{array}
[c]{l}%
\text{associative }\\
\text{cycles}%
\end{array}
$ & \multicolumn{1}{|c|}{$\Omega_{\mathcal{M}^{ass}\left(  M\right)  }$} &
\multicolumn{1}{|c|}{$\Theta_{\mathcal{M}^{ass}\left(  M\right)  }$} &
\multicolumn{1}{|c|}{$\mathcal{Y}_{\mathcal{M}^{ass}\left(  M\right)  }$%
}\\\hline
$\mathcal{M}^{bdl}\left(  M\right)  $ & $%
\begin{array}
[c]{l}%
\text{deformed DT }\\
\text{bundles}%
\end{array}
$ & \multicolumn{1}{|c|}{$\Omega_{\mathcal{M}^{bdl}\left(  M\right)  }$} &
\multicolumn{1}{|c|}{$\Theta_{\mathcal{M}^{bdl}\left(  M\right)  }$} &
\multicolumn{1}{|c|}{$\mathcal{Y}_{\mathcal{M}^{bdl}\left(  M\right)  }$%
}\\\hline
\end{tabular}

{\small Table: Structures on the moduli space of various cycles on M.}
\end{center}

\bigskip

In section \ref{1SecFibFourier} we study structures of calibrated fibrations
on $M$, they are associative $T^{3}$-fibration, coassociative $T^{4}%
$-fibration and coassociative K3-fibration. Then we propose a duality
transformation along such fibrations on $G_{2}$-manifolds analogous to the
Fourier-Mukai transformation in the geometric mirror symmetry conjecture by
Strominger, Yau and Zaslow \cite{SYZ}. Roughly speaking it should be given by
a fiberwise Fourier transformation on a coassociative $T^{4}$-fibration on
$M$. We partially verify our proposed conjecture in the flat case, in the
spirit of \cite{Le} and \cite{LYZ}.\footnote{Those authors are also familiar
with these transformations in the $G_{2}$ case.}

For instance the transformation of an ASD connection over a coassociative
torus fiber is another ASD connection over the dual torus fiber, as studied
earlier by Schenk \cite{Sc}, Braam and van Baal \cite{BB}.

We also study the Fourier transformation on an associative $T^{3}$-fibration
on $M$. The Fourier transformation on a coassociative $T^{4}$-fibration takes
cycles calibrated by $e^{\Theta}$\ (resp. $\ast e^{\Theta}$) back to
themselves. On an associative $T^{3}$-fibration, the Fourier transformation
takes cycles calibrated by $e^{\Theta}$ to those calibrated by $\ast
e^{\Theta}$.

In section \ref{Sec Spin(7)}, we study the geometry of $Spin\left(  7\right)
$-manifolds in a same manner, and therefore our discussions will be brief.

\section{$G_{2}$-manifolds and their moduli}

\subsection{Definitions}

We first review some basic structures on $G_{2}$-manifolds (see \cite{Jo} or
\cite{Sa} for more details). Any $G_{2}$-manifold $M$ admits a parallel
positive three form $\Omega_{M}$ (or simply $\Omega$) and its Ricci curvature
is zero. When
\[
M=\mathbb{R}^{7}=\operatorname{Im}\mathbb{H}\oplus\mathbb{H}=\operatorname{Im}%
\mathbb{O}%
\]
with the standard metric and orientation $dx^{123}dy^{0123}$we have,
\[
\Omega=dx^{123}-dx^{1}\left(  dy^{23}+dy^{10}\right)  -dx^{2}\left(
dy^{31}+dy^{20}\right)  -dx^{3}\left(  dy^{12}+dy^{30}\right)  .
\]

Remark: Since the Lie group $G_{2}$ preserves a vector cross product structure
on $\mathbb{R}^{7}$, any $G_{2}$-manifold inherits such a product structure
$\times$ on its tangent spaces given by,
\[
\left\langle u\times v,w\right\rangle =\Omega\left(  u,v,w\right)  .
\]
The parallel four form $\ast\Omega$ will be denoted as $\Theta$,%

\[
\Theta=dy^{0123}+dx^{23}\left(  dy^{23}+dy^{10}\right)  +dx^{31}\left(
dy^{31}+dy^{20}\right)  +dx^{12}\left(  dy^{12}+dy^{30}\right)  \text{.}%
\]

When $M$ is the product of a Calabi-Yau threefold $X$ and $S^{1}$ with its
reversed orientation, then the holonomy group equals $Hol\left(  M\right)
=SU\left(  3\right)  \subset G_{2}$. Thus $M=X\times S^{1}$ is a $G_{2}%
$-manifold. In this case we have
\begin{align*}
\Omega &  =\operatorname{Re}\Omega_{X}-\omega_{X}\wedge d\theta\\
\Theta &  =\operatorname{Im}\Omega_{X}\wedge d\theta-\omega_{X}^{2}/2\text{.}%
\end{align*}
where $\omega_{X}$ and $\Omega_{X}$ are the Ricci flat K\"{a}hler form and a
holomorphic volume form on $X$.

Using the $G_{2}$ action, we can decompose differential forms (or cohomology
classes) on $M$ into irreducible components. For example,
\begin{align*}
\Lambda^{1}  &  =\Lambda_{7}^{1}\\
\Lambda^{2}  &  =\Lambda_{7}^{2}\oplus\Lambda_{14}^{2}\\
\Lambda^{3}  &  =\Lambda_{1}^{3}\oplus\Lambda_{7}^{3}\oplus\Lambda_{27}^{3}.
\end{align*}
When $Hol\left(  M\right)  =G_{2}$, we have $H_{7}^{k}\left(  M\right)  =0$
for all $k$ and $H_{1}^{3}\left(  M\right)  $ (resp. $H_{1}^{4}\left(
M\right)  $) is generated by $\Omega$ (resp. $\Theta$). Furthermore the moduli
space of $G_{2}$-metrics is smooth with tangent space equals $H^{3}\left(
M\right)  =H_{1}^{3}\left(  M\right)  \oplus H_{27}^{3}\left(  M\right)  $. If
we normalize the volume to be one, the tangent space to this moduli space
$\mathcal{M}^{G_{2}}$ becomes $H_{27}^{3}\left(  M\right)  $.

\subsection{$G_{2}$-Analog of Yukawa couplings}

In this subsection we introduce several symmetric tensors on the moduli space
of $G_{2}$-manifolds. They are the Yukawa coupling, a metric tensor and the
prepotential. First we recall a tensor $\chi\in\Omega^{3}\left(
M,T_{M}\right)  $ which is defined by $\iota_{v}\Theta=\left\langle
\chi,v\right\rangle \in\Omega^{3}\left(  M\right)  $ for every tangent vector
$v$.

\begin{definition}
Suppose $M$ is a compact $G_{2}$-manifold, we define
\[
\mathcal{C}_{M}:%
{\textstyle\bigotimes^{3}}
H^{3}\left(  M,\mathbb{R}\right)  \rightarrow\mathbb{R}%
\]%
\[
\mathcal{C}_{M}\left(  \phi_{1},\phi_{2},\phi_{3}\right)  =\int_{M}%
\Omega\left(  \hat{\phi}_{1},\hat{\phi}_{2},\hat{\phi}_{3}\right)
\wedge\Theta,
\]
where $\hat{\phi}=\ast\left(  \phi\wedge\chi\right)  \in\Omega^{1}\left(
M,T_{M}\right)  $ and $\Omega\left(  \hat{\phi}_{1},\hat{\phi}_{2},\hat{\phi
}_{3}\right)  \in\Omega^{3}\left(  M\right)  $ is the evaluation of the three
form $\Omega$ on the vector components of $\hat{\phi}_{j}$'s.
\end{definition}

It is clear that $\mathcal{C}_{M}$ is symmetric. By using $\mathcal{C}_{M}$ we
define $\mathcal{Y}_{M}$, $\mathcal{G}_{M}$ and $\mathcal{F}_{M}$ on
$H_{27}^{3}\left(  M\right)  $ as follows: \footnote{Recall $H_{27}^{3}\left(
M\right)  =H^{3}\left(  M\right)  \cap Ker\left(  \wedge\Omega\right)  \cap
Ker\left(  \wedge\Theta\right)  .$}

\bigskip

\begin{center}%
\begin{tabular}
[c]{|l|l|}\hline
Yukawa coupling & $\mathcal{Y}_{M}\left(  \phi\right)  =\mathcal{C}_{M}\left(
\phi,\phi,\phi\right)  =C\int_{M}\Omega\wedge\ast\Omega%
\begin{array}
[c]{l}%
\,\\
\,
\end{array}
$\\\hline
Metric tensor & $\mathcal{G}_{M}\left(  \phi\right)  =\mathcal{C}_{M}\left(
\phi,\phi,\Omega\right)  =C^{\prime}\int_{M}\phi\wedge\ast\phi%
\begin{array}
[c]{l}%
\,\\
\,
\end{array}
$\\\hline
Prepotential & $\mathcal{F}_{M}=\mathcal{C}_{M}\left(  \Omega,\Omega
,\Omega\right)  .
\begin{array}
[c]{l}%
\,\\
\,
\end{array}
$\\\hline
\end{tabular}

{\small Table: Symmetric tensors on the moduli space of G}$_{2}$%
{\small -manifolds. Here }$\phi\in H_{27}^{3}\left(  M\right)  .$
\end{center}

We will use this Yukawa coupling and $p_{1}\left(  M\right)  $ in section
\ref{1SecCalibFib} to describe structures of calibrated fibrations on $M$.

\bigskip

\textbf{Reduction to CY threefolds}

When $M=X\times S^{1}$ with $Hol\left(  X\right)  =SU\left(  3\right)  $, i.e.
$X $ is a Calabi-Yau threefold, any deformation of the Einstein structure on
$M$ remains a product. There is a natural isomorphism of complex vector
spaces,
\[
H_{27}^{3}\left(  M,\mathbb{C}\right)  \cong H^{1,1}\left(  X\right)
+H^{2,1}\left(  X\right)  +H^{1,2}\left(  X\right)  \text{.}%
\]
For example if $\phi$ is a primitive form in $H^{1,1}\left(  X\right)  $
(resp. $H^{2,1}\left(  X\right)  $) then $\phi\wedge dt$ (resp. $\phi$) is in
$H_{27}^{3}\left(  M,\mathbb{C}\right)  $. Moreover $\mathcal{Y}_{M}$
restricts to the Yukawa couplings on $H^{1,1}\left(  X\right)  $,
$H^{2,1}\left(  X\right)  $ and $H^{1,2}\left(  X\right)  $ (e.g. \cite{LYZ}).

\bigskip

Remark: The cubic structure on the $G_{2}$-module $\Lambda_{27}^{3}$ is the
exceptional Jordan algebra structure on $\mathbb{R}^{27}$, whose automorphism
group determines an exceptional Lie group of type $E_{6}$. When we reduce from
$G_{2}$ to $SU\left(  3\right)  $, i.e. $M=X\times S^{1}$, the K\"{a}hler form
$\omega_{X}$ is naturally an element in $\Lambda_{27}^{3}$ which splits
$\Lambda_{27}^{3}=\mathbf{27}$ into $\mathbf{1+26}$. In this case the
automorphism group of the algebraic structure of this 26 dimensional vector
space determines an exceptional Lie group of type $F_{4}$. Roughly speaking
the Yukawa coupling on a $G_{2}$-manifold arises from a $E_{6}$ structure.
When the holonomy group reduces to $SU\left(  3\right)  $, then the $E_{6}$
structure also reduces to a $F_{4}$ structure.

\bigskip

\textbf{Including }$B$\textbf{-fields}

From physical considerations, Acharya \cite{Ac} argues that when $M$ is a
smooth $G_{2}$ manifold, then the low energy theory for M-theory on
$M\times\mathbb{R}^{3,1}$ is a $N=1$ supergravity theory with $b_{2}\left(
M\right)  $ Abelian vector multiplets plus $b_{3}\left(  M\right)  $ neutral
chiral multiplets. We also include one scalar field and consider the enlarged
moduli space $\mathcal{M}^{G_{2}+B}=\mathcal{M}^{G_{2}}\times H^{2}\left(
M,U\left(  1\right)  \right)  \times H^{0}\left(  M,U\left(  1\right)
\right)  $. Note $b_{0}\left(  M\right)  =1$.

We are going to construct the corresponding Yukawa coupling for this moduli
space. We define (i) a symmetric cubic tensor $\mathcal{C}^{\prime}$ on
$H^{2}\left(  M,\mathbb{R}\right)  $,
\begin{align*}
\mathcal{C}_{M}^{\prime}  &  :%
{\textstyle\bigotimes^{3}}
H^{2}\left(  M\right)  \rightarrow\mathbb{R}\\
\mathcal{C}_{M}^{\prime}\left(  \beta_{1},\beta_{2},\beta_{3}\right)   &
=\int_{M}\Omega\left(  \hat{\beta}_{1},\hat{\beta}_{2},\hat{\beta}_{3}\right)
\wedge\Theta,
\end{align*}
where $\hat{\beta}=\iota_{\beta}\chi\in\Omega^{1}\left(  M,T_{M}\right)  $ and
(ii) a symmetric bilinear tensor $\mathcal{Q}_{M}$ on $H^{2}\left(
M,\mathbb{R}\right)  ,$%
\begin{align*}
\mathcal{Q}_{M}  &  :%
{\textstyle\bigotimes^{2}}
H^{2}\left(  M\right)  \rightarrow\mathbb{R}\\
\mathcal{Q}_{M}\left(  \beta_{1},\beta_{2}\right)   &  =\int_{M}\beta
_{1}\wedge\beta_{2}\wedge\Omega.
\end{align*}

Now the Yukawa coupling on the enlarged moduli space $\mathcal{M}^{G_{2}%
}\times H^{2}\left(  M,U\left(  1\right)  \right)  \times H^{0}\left(
M,U\left(  1\right)  \right)  $ of $G_{2}$-manifolds (with $Hol=G_{2}$) is
defined as a combination of $\mathcal{Y}_{M}$, $\mathcal{C}_{M}^{\prime}$ and
$\mathcal{Q}_{M}$ as follow,
\[
\mathcal{\tilde{Y}}_{M}:%
{\textstyle\bigotimes^{3}}
\left(  H_{27}^{3}\left(  M\right)  +H_{14}^{2}\left(  M\right)  +H^{0}\left(
M\right)  \right)  \rightarrow\mathbb{R}\text{.}%
\]

When $M=X\times S^{1}$ for a Calabi-Yau threefold $X$ we have $H_{14}%
^{2}\left(  M\right)  =H_{prim}^{1,1}\left(  X,\mathbb{R}\right)  $ which
corresponds to the decomposition $\mathbf{g}_{2}=\mathbf{su}\left(  3\right)
+\mathbb{C}^{3}+\mathbb{C}^{3\ast}$ because of $H^{2,0}\left(  X\right)
=H^{0,2}\left(  X\right)  =0$ for a manifold with $SU\left(  3\right)  $
holonomy. Therefore
\[
H_{14}^{2}\left(  M,U\left(  1\right)  \right)  +H^{0}\left(  M,U\left(
1\right)  \right)  \cong H^{2}\left(  X,U\left(  1\right)  \right)  \text{,}%
\]
and an element in it is usually called a \textit{B-field} on the Calabi-Yau
manifold $X$. Our Yukawa coupling $\mathcal{Y}^{\prime}$ thus reduces to the
usual one on Calabi-Yau manifolds coupled with B-fields.

\bigskip

Remark: If we do not fix the volume of $M$, then the moduli space of $G_{2}%
$-metrics with $B$-fields is locally isomorphic to $H^{\leq3}\left(
M,\mathbb{R}\right)  $. In M-theory (see e.g. \cite{AW}) it is desirable to
complexify this moduli space and this implies that the resulting space will be
locally isomorphic to the total cohomology group $H^{\ast}\left(
M,\mathbb{R}\right)  $ of $M$.

\section{Supersymmetric cycles}

\subsection{\label{1SecDTbdl}Deformed Donaldson-Thomas bundles}

In \cite{DT}, Donaldson and Thomas study the following Yang-Mills equation for
Hermitian connections $D_{E}$ on a complex vector bundle $E$ over $M$,
\[
F_{E}\wedge\Theta=0\in\Omega^{6}\left(  M,ad\left(  E\right)  \right)
\text{.}%
\]
This is the Euler-Lagrange equation for the following Chern-Simons type
functional,
\begin{align*}
\mathcal{A}\left(  E\right)   &  \rightarrow\mathbb{R}\\
D_{E}  &  \rightarrow\int_{M}CS_{3}\left(  D_{0},D_{E}\right)  \wedge
\Theta\text{,}%
\end{align*}
where $D_{0}$ is any fixed background connection. This functional is
equivalent to the following functional
\[
\int_{M\times\left[  0,1\right]  }Tre^{\tilde{F}}\wedge\Theta,
\]
where $\tilde{F}$ is the curvature of a connection on $M\times\left[
0,1\right]  $ formed by the affine path of connections on $E$ joining $D_{0}$
and $D_{E}$.

\bigskip

On a $G_{2}$-manifold we introduce the following \textit{deformed
Donaldson-Thomas equation} which should have the corresponding effect of
preserving supersymmetry in M-theory \cite{MMMS},%

\[
F_{E}\wedge\Theta+F_{E}^{3}/6=0\in\Omega^{6}\left(  M,ad\left(  E\right)
\right)  .
\]
Equivalently this equals,
\[
\left[  e^{\Theta+F_{E}}\right]  ^{\left[  6\right]  }=0\in\Omega^{6}\left(
M,ad\left(  E\right)  \right)  .
\]
It is the Euler-Lagrangian equation for the following Chern-Simons type
functional,
\[
D_{E}\rightarrow\int_{M\times\left[  0,1\right]  }Tr\left[  e^{\Theta
+\tilde{F}}\right]  =\int_{M}CS\left(  D_{0},D_{E}\right)  e^{\Theta}.
\]

The mirror transformation along a coassociative $T^{4}$-fibration on a flat
$G_{2}$-manifold, which we will discuss in section \ref{Sec Trans coa}, takes
an associative section to a deformed Donaldson-Thomas bundles on the mirror manifold.

\bigskip

\textbf{Geometric structures on }$\mathcal{M}^{bdl}\left(  M\right)  $

We define natural differential forms and a symmetric cubic tensor on the
moduli space of deformed DT bundles $\mathcal{M}^{bdl}\left(  M\right)  $ as follows,

\begin{center}
\bigskip%
\begin{tabular}
[c]{|l|l|}\hline
3-form: & $\Omega_{\mathcal{M}^{bdl}\left(  M\right)  }=\int_{M}Tr\left[
\alpha\wedge\beta\wedge\gamma\right]  _{skew}e^{\Theta+F_{E}}
\begin{array}
[c]{l}%
\,\\
\,
\end{array}
$\\\hline
4-form: & $\Theta_{\mathcal{M}^{bdl}\left(  M\right)  }=\int_{M}Tr\left[
\alpha\wedge\beta\wedge\gamma\wedge\delta\right]  _{skew}\ast e^{\Theta
+F_{E}}
\begin{array}
[c]{l}%
\,\\
\,
\end{array}
$\\\hline
Cubic: & $\mathcal{Y}_{\mathcal{M}^{bdl}\left(  M\right)  }=\int
_{M}\left\langle \alpha,\left[  \beta,\gamma\right]  \right\rangle
_{ad}e^{\Theta+F_{E}}.
\begin{array}
[c]{l}%
\,\\
\,
\end{array}
$\\\hline
\end{tabular}

{\small Table: Natural geometric structures on the moduli space of deformed }

{\small DT bundles on M. Here }$\alpha,\beta,\gamma,\delta\in\Omega^{1}\left(
M,ad\left(  E\right)  \right)  .$
\end{center}

\bigskip

Here $\left\langle \cdot,\cdot\right\rangle _{ad}$ is the Killing form on
$\mathbf{g}$, the Lie algebra of the gauge group. In particular the symmetric
cubic tensor $\mathcal{Y}_{\mathcal{M}^{bdl}\left(  M\right)  }$ on
$\mathcal{M}^{bdl}\left(  M\right)  $ is trivial when $G$ is Abelian. Note
$e^{\Theta+F_{E}}$ in the cubic form works as $\left[  e^{\Theta+F_{E}%
}\right]  ^{\left[  6\right]  }$.

\bigskip

Remark: When $M$ is a flat torus $T=\mathbb{R}^{7}/\Lambda$, the moduli space
$\mathcal{M}^{bdl}\left(  M\right)  $ of flat $U\left(  1\right)  $-bundles is
canonically isomorphic to the dual flat torus $T^{\ast}=\mathbb{R}^{7\ast
}/\Lambda^{\ast}$. Under this natural identification, we have
\begin{align*}
\Omega_{\mathcal{M}^{bdl}\left(  T\right)  }  &  =\Omega_{T^{\ast}},\\
\Theta_{\mathcal{M}^{bdl}\left(  T\right)  }  &  =\Theta_{T^{\ast}}.
\end{align*}
Moreover this transformation from $T$ to $\mathcal{M}^{bdl}=T^{\ast}$ is
involutive, i.e. $\mathcal{M}^{bdl}\left(  T^{\ast}\right)  =T$.

\bigskip

Remark: The tangent bundle of $M$ is a $G_{2}$-bundle, moreover the curvature
tensor of the Levi-Civita connection on $M$ satisfies the DT-equation,
\[
F\wedge\Theta=0\text{.}%
\]
This equation is equivalent to $F\in\Omega_{14}^{2}\left(  M,ad\left(
T_{M}\right)  \right)  $ and it follows from the fact that $F\in\Omega
^{2}\left(  M,\mathbf{g}_{2}\left(  T_{M}\right)  \right)  $ and the torsion
freeness of the connection.

Infinitesimal deformations of $T_{M}$ as a DT-bundle are parametrized by the
first cohomology group $H^{1}\left(  M,\mathbf{g}_{2}\left(  T_{M}\right)
\right)  $ of the following elliptic complex \cite{DT}:
\[
0\rightarrow\Omega^{0}\left(  M,\mathbf{g}_{2}\left(  T_{M}\right)  \right)
\overset{D_{E}}{\rightarrow}\Omega^{1}\left(  M,\mathbf{g}_{2}\left(
T_{M}\right)  \right)  \overset{\Theta D_{E}}{\rightarrow}\Omega^{6}\left(
M,\mathbf{g}_{2}\left(  T_{M}\right)  \right)  \overset{D_{E}}{\rightarrow
}\Omega^{7}\left(  M,\mathbf{g}_{2}\left(  T_{M}\right)  \right)
\rightarrow0.
\]

When $M=X\times S^{1}$ with $X$ a Calabi-Yau threefold with $SU\left(
3\right)  $ holonomy, then there is a natural isomorphism,
\[
H^{1}\left(  M,\mathbf{g}_{2}\left(  T_{M}\right)  \right)  \cong
H^{1,1}\left(  X\right)  +H^{2,1}\left(  X\right)  +H^{1}\left(
X,End_{0}\left(  T_{X}\right)  \right)  .
\]
Furthermore the above cubic tensor $\mathcal{Y}_{\mathcal{M}^{bdl}\left(
M\right)  }$ restricts to $\mathcal{Y}_{A}$ on $H^{1,1}\left(  X\right)  $,
$\mathcal{Y}_{B}$ on $H^{2,1}\left(  X\right)  $ and a similar Yukawa coupling
$\mathcal{Y}_{C}$ on $H^{1}\left(  X,End_{0}\left(  T_{X}\right)  \right)  $.
All three seemingly unrelated Yukawa couplings on a Calabi-Yau threefold
combine in a natural way when we study the $G_{2}$-manifold $X\times S^{1}$!

\subsection{\label{1SecAssCycle}\label{Sec Ass cycle}Associative cycles}

\textbf{Associative submanifolds}

There are two types of nontrivial calibrations on a $G_{2}$-manifold $M$ as
studied by Harvey and Lawson (\cite{HL}), they are (i) associative
submanifolds and (ii) coassociative submanifolds. They are always absolute
minimal submanifolds in $M$.

A three dimensional submanifold $A$ in $M$ is called an \textit{associative}
submanifold if the restriction of $\Omega$ to $A$ equals the volume form on
$A$ of the induced metric. It has the following two equivalent
characterizations: (1) the restriction of $\chi\in\Omega^{3}\left(
M,T_{M}\right)  $ to $A$ is zero;
\[
\chi|_{A}=0\in\Omega^{3}\left(  A,T_{M}|_{A}\right)  \text{.}%
\]
(2) the vector cross product on $M$ preserves $T_{A}$, i.e. $u,v\in T_{A}$
implies that $u\times v\in T_{A}$. As a corollary, two associative
submanifolds can not intersect along a two dimensional subspace.

McLean studies the deformation theory of associative submanifolds in \cite{Mc}
and identifies their infinitesimal deformations as certain twisted harmonic
spinors. Unlike coassociative submanifolds, deformations of an associative
submanifold can be \textit{obstructed}. For example, if we take any smooth
isolated rational curve $C$ with normal bundle $O\oplus O\left(  -2\right)  $
in a Calabi-Yau threefold $X$, then $C\times S^{1}$ is an associative
submanifold in $M=X\times S^{1}$ with obstructed deformations. That is the
moduli space of associative submanifolds, denoted $\mathcal{B}^{ass}\left(
M\right)  $, can be singular.

\bigskip

\textbf{Associative cycles}

From a physical perspective, one would couple associative submanifolds with
gauge fields. This is analogous to coupling of gauge fields with special
Lagrangian submanifolds of a Calabi-Yau manifold. From a mathematical
perspective, as we will see, the moduli space of associative submanifolds
coupled with gauge fields has very rich geometric structures.

An \textit{associative cycle }is defined to be any pair $\left(
A,D_{E}\right)  $ with $A$ an associative submanifold in $M$ and $D_{E}$ a
unitary flat connection on a bundle $E$ over $A$. The flatness equation,
$F_{E}=0$ for connections over a three manifold is the Euler-Lagrange equation
of the standard Chern-Simons functional. Analogously the associativity
condition for a pair $\left(  A,D_{E}\right)  $ is the Euler-Lagrange equation
of the following functional:
\begin{align*}
Map\left(  A,M\right)  \times\mathcal{A}\left(  E\right)   &  \rightarrow
\mathbb{R}/\mathbb{Z}\\
CS\left(  A,D_{E}\right)   &  =\int_{A\times\left[  0,1\right]  }Tr\left[
e^{\Theta+\tilde{F}}\right]  .
\end{align*}
The first term in the above functional is a direct analog of a functional used
by Thomas \cite{Th} for special Lagrangian submanifolds in a Calabi-Yau threefold.

Next we are going to study the moduli space of associative cycles, namely the
critical set of $CS$, and we denote it as $\mathcal{M}^{ass}\left(  M\right)
$. This space has richer structure than the moduli space of associative
submanifolds $\mathcal{B}^{ass}\left(  M\right)  $.

The tangent space of $\mathcal{M}^{ass}\left(  M\right)  $ at an associative
cycle $\left(  A,D_{E}\right)  $ can be identified as (see \cite{Mc}),
\[
T_{\left(  A,D_{E}\right)  }\mathcal{M}^{ass}=Ker\mathbf{D}\oplus H^{1}\left(
A,ad\left(  E\right)  \right)  .
\]
We define a three form $\Omega_{\mathcal{M}^{ass}\left(  M\right)  }$, four
form $\Theta_{\mathcal{M}^{ass}\left(  M\right)  }$ and a symmetric cubic
tensor $\mathcal{Y}_{\mathcal{M}^{ass}\left(  M\right)  }$ on $\mathcal{M}%
^{ass}\left(  M\right)  $ as follows:

\begin{center}%
\begin{tabular}
[c]{|l|l|}\hline
3-form: & $\Omega_{\mathcal{M}^{ass}\left(  M\right)  }=\left\{
\begin{array}
[c]{l}%
\,\int_{A}Tr\alpha\wedge\beta\wedge\gamma\medskip\,\\
-\,\int_{A}\left\langle \alpha\cdot\overline{\phi},\overline{\eta
}\right\rangle \Omega
\end{array}
\right.
\begin{array}
[c]{l}%
\,\\
\\
\,
\end{array}
$\\\hline
4-form: & $\Theta_{\mathcal{M}^{ass}\left(  M\right)  }\,=\left\{
\begin{array}
[c]{l}%
\,\int_{A}\left\langle \overline{\phi},\ast\left(  \alpha\wedge\beta\right)
\cdot\overline{\eta}\right\rangle \Omega\medskip\,\\
\int_{A}\det\left(  \phi,\eta,\xi,\zeta\right)  \Omega\medskip
\end{array}
\right.
\begin{array}
[c]{l}%
\,\\
\\
\,
\end{array}
$\\\hline
Cubic: & $\mathcal{Y}_{\mathcal{M}^{ass}\left(  M\right)  }=\int
_{A}\left\langle \left[  \alpha,\beta\right]  ,\gamma\right\rangle \Omega$.$%
\begin{array}
[c]{l}%
\,\\
\,
\end{array}
$\\\hline
\end{tabular}

{\small Table: Natural geometric structures on the moduli space of associative
cycles}

{\small on M. Here }$\alpha,\beta,\gamma,\delta\in\Omega^{1}\left(
A,ad\left(  E\right)  \right)  $ {\small and }$\phi,\eta,\xi,\zeta\in
KerD${\small .}
\end{center}

In the definition of the four form, we use the fact that the twisted spinor
bundle $S$ is a complex rank two bundle and therefore, as a real vector bundle
$S_{\mathbb{R}}$, there is a natural trivialization of $\Lambda^{4}%
S_{\mathbb{R}}\cong\mathbb{R}$. The action of $\operatorname{Im}\mathbb{H}$ on
$\mathbb{H}$ has used to define $3$ and $4$-form.

\bigskip

\textbf{An example}

In this example, the above three form $\Omega_{\mathcal{M}^{ass}\left(
M\right)  }$ defines a $G_{2}$ structure on the moduli space $\mathcal{M}%
^{ass}\left(  M\right)  $. We consider a flat example $M=T^{3}\times T^{4}$
with $T^{3}=\mathbb{R}^{3}/\Lambda_{3}$ and $T^{4}=\mathbb{R}^{4}/\Lambda_{4}%
$. The projection to the second factor $M\rightarrow B=T^{4}$ is an
associative fibration. The moduli space $\mathcal{M}^{ass}\left(  M\right)  $
of $\left(  A,D_{E}\right)  $ with $A$ a fiber and $D_{E}$ a flat $U\left(
1\right)  $-connection over $A$, can be identified with $T^{3\ast}\times
T^{4}$ where $T^{3\ast}=\mathbb{R}^{3\ast}/\Lambda_{3}^{\ast}$ is the dual
torus to $T^{3}$. Furthermore $\Omega_{\mathcal{M}^{ass}\left(  M\right)  }$
and $\Theta_{\mathcal{M}^{ass}\left(  M\right)  }$ are precisely the natural
calibration three and four forms on the flat $G_{2}$-manifold $T^{3\ast}\times
T^{4}$. This is a simple but fun exercise.

In particular, $\chi_{\mathcal{M}^{ass}\left(  T^{3}\times T^{4}\right)
}=\chi_{\mathcal{M}^{ass}\left(  T^{3\ast}\times T^{4}\right)  }$ and
therefore, $T^{3\ast}\times T^{4}\rightarrow T^{4}$ defines an associative
fibration structure on $\mathcal{M}^{ass}\left(  T^{3}\times T^{4}\right)
\rightarrow\mathcal{B}^{ass}\left(  T^{3}\times T^{4}\right)  $ with a
coassociative section, namely, $\Omega_{\mathcal{M}^{ass}\left(  T^{3}\times
T^{4}\right)  }$ (resp. $\chi_{\mathcal{M}^{ass}\left(  T^{3}\times
T^{4}\right)  }$) vanishes on the section (resp. fibers) of this fibration.

In the next paragraph, we will explain that such structures exist on every
$\mathcal{M}^{ass}\left(  M\right)  $. However, when $M$ is not flat,
$\Omega_{\mathcal{M}^{ass}\left(  M\right)  }$ and $\Theta_{\mathcal{M}%
^{ass}\left(  M\right)  }$ are not parallel forms, thus they do not define
$G_{2}$-structure on $\mathcal{M}^{ass}\left(  M\right)  $ and we can only
discuss their (co-)associativity in terms of vanishing of appropriate tensors.

\bigskip

\textbf{Moduli space of associative submanifolds} \textbf{is coassociative}

Every associative submanifold $A$ defines an associative pair $\left(
A,D_{E}\right)  $ where $D_{E}$ is simply the trivial connection $d$. Thus we
have an embedding of $\mathcal{B}^{ass}\left(  M\right)  $ inside
$\mathcal{M}^{ass}\left(  M\right)  $. It is not difficult to see that the
restriction of $\Omega_{\mathcal{M}^{ass}\left(  M\right)  }$ to it is zero.
As we will see in the next subsection, this property characterizes a
coassociative submanifold is $\Omega_{\mathcal{M}^{ass}\left(  M\right)  }$
defines a $G_{2} $ structure on $\mathcal{M}^{ass}\left(  M\right)  $.
Furthermore we have a natural fibration structure, $\mathcal{M}^{ass}\left(
M\right)  \rightarrow\mathcal{B}^{ass}\left(  M\right)  $ by forgetting the
connections. Fibers are associative submanifolds in the sense that the
restriction $\chi_{\mathcal{M}^{ass}\left(  M\right)  }$ to them are zero.
Here $\chi_{\mathcal{M}^{ass}\left(  M\right)  }$ is defined the same way as
$\chi$ for $M$. The reason is fiber directions for $\mathcal{M}^{ass}\left(
M\right)  \rightarrow\mathcal{B}^{ass}\left(  M\right)  $ corresponds to
$H^{1}\left(  A,ad\left(  E\right)  \right)  $ and each component in
$\Theta_{\mathcal{M}^{ass}\left(  M\right)  }$ involves at least 2
$Ker\mathbf{D}$ components in its variables. Therefore $\chi_{\mathcal{M}%
^{ass}}$ must vanish on fibers.

\subsection{\label{Sec Coa cycle}Coassociative cycles}

\textbf{Coassociative submanifolds}

A four dimensional submanifold $C$ in $M$ is called a \textit{coassociative
submanifold} if it is calibrated by $\Theta$, i.e. $\Theta|_{C}=vol_{C}$. It
can be characterized by
\[
\Omega|_{C}=0\in\Omega^{3}\left(  C,\mathbb{R}\right)  .
\]
Another characterization is the vector cross product on $M$ preserves
$N_{C/M}\subset TM$. As a corollary of this, any two coassociative
submanifolds in $M$ can not intersect along a three dimension subspace.

The normal bundle of $C$ in $M$ can be identified with the bundle of
anti-self-dual two forms on $C$,
\[
N_{C/M}=\Lambda_{+}^{2}\left(  C\right)  .
\]
Furthermore the space of infinitesimal deformations of $C$ inside $M$ as a
coassociative submanifold is isomorphic to the space of self-dual harmonic two
forms on $C$, $H_{+}^{2}\left(  C,\mathbb{R}\right)  $. In fact the
deformation theory of coassociative submanifold is unobstructed \cite{Mc} and
therefore their moduli space, denoted $\mathcal{B}^{coa}\left(  M\right)  $,
is always smooth.

\bigskip

\textbf{Coassociative cycles}

Similar to the associative cases, we will also couple coassociative
submanifolds with certain gauge fields. A \textit{coassociative cycle} is a
pair $\left(  C,D_{E}\right)  $ with $C$ a coassociative submanifold in $M$
and $E$ an anti-self-dual (ASD) connection on $C$ with respect to the induced
metric. Recall that a connection is called ASD if its curvature $2$-form
$F_{E}$ lies in $\Omega_{-}^{2}\left(  C,ad\left(  E\right)  \right)  $, i.e.
$F_{E}+$ $\ast F_{E}=0$. The tangent space of their moduli space
$\mathcal{M}^{coa}\left(  M\right)  $ can be identified as $H_{+}^{2}\left(
C\right)  +H^{1}\left(  C,ad\left(  E\right)  \right)  $ \cite{Mc}. We define
a three form $\Omega_{\mathcal{M}^{coa}\left(  M\right)  }$, a four form
$\Theta_{\mathcal{M}^{coa}\left(  M\right)  }$ and a symmetric cubic tensor
$\mathcal{Y}_{\mathcal{M}^{coa}\left(  M\right)  }$ on $\mathcal{M}%
^{coa}\left(  M\right)  $ as follows:

\bigskip

\begin{center}%
\begin{tabular}
[c]{|l|l|}\hline
3-form: & $\Omega_{\mathcal{M}^{coa}\left(  M\right)  }=\left\{
\begin{array}
[c]{l}%
\int_{C}\left[  \phi,\eta\right]  \wedge\xi\medskip\,\\
-\int_{C}Tr\left(  \phi\wedge\alpha\wedge\beta\right)
\end{array}
\right.
\begin{array}
[c]{l}%
\,\\
\\
\,
\end{array}
$\\\hline
4-form: & $\Theta_{\mathcal{M}^{coa}\left(  M\right)  }\,=\left\{
\begin{array}
[c]{l}%
-\int_{C}Tr\left(  \alpha\wedge\beta\wedge\gamma\wedge\delta\right)
_{skew}\medskip\,\\
\int_{C}Tr\left(  \left[  \phi,\eta\right]  \wedge\alpha\wedge\beta\right)
\end{array}
\right.
\begin{array}
[c]{l}%
\,\\
\\
\,
\end{array}
$\\\hline
Cubic: & $\mathcal{Y}_{\mathcal{M}^{coa}\left(  M\right)  }=\int_{A}%
Tr\phi\left[  \alpha,\beta\right]  _{ad\left(  E\right)  }$.$%
\begin{array}
[c]{l}%
\,\\
\,
\end{array}
$\\\hline
\end{tabular}

{\small Table: Natural geometric structures on the moduli space of
coassociative cycles on M. Here }$\alpha,\beta,\gamma,\delta\in\Omega
^{1}\left(  C,ad\left(  E\right)  \right)  $ {\small and }$\phi,\eta,\xi\in
H_{+}^{2}\left(  C\right)  ${\small .}
\end{center}

Here, Lie algebra structure on space of $2$-forms is given by $\wedge
^{2}\mathbb{R}^{4}\simeq so\left(  \mathbb{R}^{4}\right)  $.

\bigskip

\textbf{An example}

In this example, the above three form $\Omega_{\mathcal{M}^{coa}\left(
M\right)  }$ defines a $G_{2}$ structure on the moduli space $\mathcal{M}%
^{coa}\left(  M\right)  $. We consider a flat example $M=T^{3}\times T^{4}$
with $T^{3}=\mathbb{R}^{3}/\Lambda_{3}$ and $T^{4}=\mathbb{R}^{4}/\Lambda_{4}%
$. The projection to the first factor $M\rightarrow B=T^{3}$ is a
coassociative fibration. The moduli space $\mathcal{M}^{coa}$ of $\left(
C,D_{E}\right)  $ with $C$ a fiber and $D_{E}$ a flat $U\left(  1\right)
$-connection over $C$, can be identified with $T^{3}\times T^{4\ast}$ where
$T^{4\ast}=\mathbb{R}^{4\ast}/\Lambda_{4}^{\ast}$ is the dual torus to $T^{4}%
$. Moreover $\Omega_{\mathcal{M}^{coa}}$ and $\Theta_{\mathcal{M}^{coa}}$ are
precisely the natural calibration three and four form on the flat $G_{2}%
$-manifold $T^{3}\times T^{4\ast}$. This is another simple but fun exercise.

\bigskip

\textbf{Moduli space of coassociative submanifolds} \textbf{is associative}

The fibration $\mathcal{M}^{coa}\left(  M\right)  \rightarrow\mathcal{B}%
^{coa}\left(  M\right)  $ is a 'coassociative fibration' and the section
embeds $\mathcal{B}^{coa}\left(  M\right)  $ inside $\mathcal{M}^{coa}\left(
M\right)  $ behaves as an 'associative' submanifold in an appropriate sense.
The arguments are parallel to the case for $\mathcal{B}^{ass}\left(  M\right)
$ in\textbf{\ }$\mathcal{M}^{ass}\left(  M\right)  $ as we discussed in
section \ref{Sec Ass cycle} and hence omitted.

\subsection{Reductions to Calabi-Yau threefolds}

In this section we assume $M=X\times S^{1}$ with $X$ a Calabi-Yau threefold
and
\[
\Omega=\operatorname{Re}\Omega_{X}-\omega\wedge dt,
\]
where $\Omega_{X}$ and $\omega_{X}$ are the holomorphic three form and the
K\"{a}hler form on $X$ respectively. We will see that every moduli space of
cycles (with $U\left(  1\right)  $-connections) on $M$ and its natural
differential forms inherit similar decompositions. That is,

\begin{center}%
\begin{tabular}
[c]{|l|}\hline
$\mathcal{M}\left(  M\right)  =\mathcal{M}\left(  X\right)  \times S^{1},
\begin{array}
[c]{l}%
\,\\
\,
\end{array}
$\\\hline
$\Omega_{\mathcal{M}\left(  M\right)  }=\operatorname{Re}\Omega_{\mathcal{M}%
\left(  X\right)  }-\omega_{\mathcal{M}\left(  X\right)  }\wedge dt.
\begin{array}
[c]{l}%
\,\\
\,
\end{array}
$\\\hline
\end{tabular}

\end{center}

\label{Decomposition}The manifold $M$ itself can be regarded as the moduli
space of zero dimensional cycles in $M$.

Remark : The reduction from a $G_{2}$-manifold $M$ to a Calabi-Yau threefold
$X$ has a physical significance. Namely we reduce the studies of M-theory on
an eleven dimensional space time $\mathbb{R}^{3.1}\times$ $M$ to the string
theory on ten dimensional space time $\mathbb{R}^{3.1}\times$ $X$.

\bigskip

\begin{center}
{\large On the moduli of deformed DT bundles}
\end{center}

The \textit{deformed Hermitian Yang-Mills }connection on a bundle $E$ on $X$
is introduced in \cite{MMMS} as a supersymmetric cycle. Its curvature tensor
satisfies the following equation,
\begin{align*}
F\wedge\omega^{2}  &  =F^{3}/3,\\
F^{0,2}  &  =0\text{.}%
\end{align*}
Their moduli space $\mathcal{M}^{bdl}\left(  X\right)  $ has a natural
holomorphic three form $\Omega_{\mathcal{M}^{bdl}\left(  X\right)  }$ and an
almost symplectic form $\omega_{\mathcal{M}^{bdl}\left(  X\right)  }$,
\begin{align*}
\Omega_{\mathcal{M}^{bdl}\left(  X\right)  }\left(  \alpha,\beta
,\gamma\right)   &  =\int_{X}\alpha\wedge\beta\wedge\gamma\wedge\Omega_{X}\\
\omega_{\mathcal{M}^{bdl}\left(  X\right)  }\left(  \alpha,\beta\right)   &
=\int_{X}\alpha\wedge\bar{\beta}\wedge e^{\omega+F_{E}}\text{.}%
\end{align*}

Here $\alpha,\beta,\gamma\in H^{1}\left(  X,O_{X}\right)  =H^{0,1}\left(
X\right)  ,$ the tangent space of $\mathcal{M}^{bdl}\left(  M\right)  $ at
$D_{E}$.

\bigskip

The pullback of any deformed Hermitian Yang-Mills connection $D_{E}$ on $X$ is
a deformed DT connection on $M$. This is because the condition $F^{0,2}=0$ for
the holomorphicity of $E$ on a threefold is equivalent to $F\wedge
\operatorname{Re}\Omega_{X}=0$ in $\Omega^{5}\left(  X\right)  $. We can
obtain other elements in $\mathcal{M}^{bdl}\left(  M\right)  $ by tensoring
with flat $U\left(  1\right)  $-connection on $S^{1}$. The moduli space of
flat connections on $S^{1}$ is the dual circle, also denoted as $S^{1}$. It is
not difficult to verify the above mentioned decompositions for $\mathcal{M}%
^{bdl}\left(  M\right)  $ and $\Omega_{\mathcal{M}^{bdl}\left(  M\right)  }$
in this case.

\begin{center}
{\large On the moduli of associative cycles}
\end{center}

Suppose that an associative submanifold $A$ in $M=X\times S^{1}$ is of product
types, then

\begin{center}%
\begin{tabular}
[c]{|l|l|}\hline
Case (1): $A=\Sigma\times S^{1}$ & $\Sigma$ holomorphic curve in $X$.$\
\begin{array}
[c]{l}%
\,\\
\,
\end{array}
$\\\hline
Case (2): $A=L\times\left\{  t\right\}  $ & $L$ special Lagrangian with phase
= 0.$%
\begin{array}
[c]{l}%
\,\\
\,
\end{array}
$\\\hline
\end{tabular}

\end{center}

Case (1): $A=\Sigma\times S^{1}$. We denote the moduli space of holomorphic
curves $\Sigma$ in $X$ as $\mathcal{B}^{cx}\left(  X\right)  $. Similarly we
consider the moduli space of pairs $\left(  \Sigma,D_{E}\right)  $ where
$D_{E} $ is a flat $U\left(  1\right)  $-connection on $\Sigma$, and denote it
as $\mathcal{M}^{cx}\left(  X\right)  $. Gopakumar and Vafa conjecture that
the cohomology group of $\mathcal{M}^{cx}\left(  X\right)  $ admits an
$\mathbf{sl}\left(  2\right)  \times\mathbf{sl}\left(  2\right)  $ action
whose multiplicities determine Gromov-Witten invariants of $X$ of every genus
for the class $\left[  \Sigma\right]  $. The tangent space of $\mathcal{M}%
^{cx}\left(  X\right)  $ at $\left(  \Sigma,D_{E}\right)  $ equals
$H^{0}\left(  \Sigma,N_{\Sigma/X}\right)  \oplus H^{1}\left(  \Sigma
,O_{\Sigma}\right)  $. It carries natural holomorphic three form
$\Omega_{\mathcal{M}^{cx}\left(  X\right)  }$ and symplectic form
$\omega_{\mathcal{M}^{cx}\left(  X\right)  }$ as follow,

\begin{center}%
\begin{tabular}
[c]{|l|l|}\hline
holomorphic 3 form: & $\Omega_{\mathcal{M}^{cx}\left(  X\right)  }%
=\int_{\Sigma}\widetilde{\left(  \phi\wedge\eta\right)  }\wedge\alpha,
\begin{array}
[c]{l}%
\,\\
\,
\end{array}
$\\\hline
symplectic form: & $\omega_{\mathcal{M}^{cx}\left(  X\right)  }=\int_{\Sigma
}\alpha\wedge\bar{\beta}\pm\left\langle \phi,\eta\right\rangle \omega_{X}.
\begin{array}
[c]{l}%
\,\\
\,
\end{array}
$\\\hline
\end{tabular}

\end{center}

Here $\alpha,\beta\in H^{1}\left(  \Sigma,O_{\Sigma}\right)  $ and $\phi
,\eta\in H^{0}\left(  \Sigma,N_{\Sigma/X}\right)  $ and $\widetilde{\left(
\phi\wedge\eta\right)  }\in T_{\Sigma}^{\ast}$ is the image $\phi\wedge\eta
\in\Lambda^{2}N_{\Sigma/X}$ under the natural identification $\Lambda
^{2}N_{\Sigma/X}\cong T_{\Sigma}^{\ast}$ because of $\Lambda^{3}T_{X}\cong
O_{X}$.

Note that any deformation of $A=\Sigma\times S^{1}$ as an associative
submanifold in $M$ must still be of the same form, i.e. $\mathcal{B}%
^{ass}\left(  M\right)  =\mathcal{B}^{cx}\left(  X\right)  $. However there
are more flat $U\left(  1\right)  $ bundles on $A$ then on $\Sigma$, i.e.
$H^{1}\left(  A,\mathbb{R}/\mathbb{Z}\right)  \cong H^{1}\left(
\Sigma,\mathbb{R}/\mathbb{Z}\right)  \times S^{1}$. We have the decompositions
for $\mathcal{M}^{coa}\left(  M\right)  $ and $\Omega_{\mathcal{M}%
^{coa}\left(  M\right)  }$.

\bigskip

Case (2): $A=L\times\left\{  t\right\}  $ with $L$ a special Lagrangian
submanifold of zero phase in $X$, i.e. calibrated by $\operatorname{Re}\Omega$.

We denote the moduli space of $L$ (resp. $\left(  L,D_{E}\right)  $ with
$D_{E} $ a flat $U\left(  1\right)  $-connection on $L$) in $X$ as
$\mathcal{B}^{cx}\left(  X\right)  $ (resp. $\mathcal{M}^{SL}\left(  X\right)
$). McLean \cite{Mc} introduces a three form on $\mathcal{B}^{cx}\left(
X\right)  $. On $\mathcal{M}^{SL}\left(  X\right)  $, there are a holomorphic
three form $\Omega_{\mathcal{M}^{SL}\left(  X\right)  }$ and a symplectic form
$\omega_{\mathcal{M}^{SL}\left(  X\right)  }$ as follow

\begin{center}%
\begin{tabular}
[c]{|l|l|}\hline
holomorphic 3 form: & $\Omega_{\mathcal{M}^{SL}\left(  X\right)  }=\ \int
_{L}\alpha\wedge\beta\wedge\gamma,
\begin{array}
[c]{l}%
\,\\
\,
\end{array}
$\\\hline
symplectic form: & $\omega_{\mathcal{M}^{SL}\left(  X\right)  }=\ \int
_{L}\left\langle \alpha,J\beta\right\rangle \operatorname{Re}\Omega_{X}.
\begin{array}
[c]{l}%
\,\\
\,
\end{array}
$\\\hline
\end{tabular}

\end{center}

Here $\alpha,\beta,\gamma\in H^{1}\left(  L,\mathbb{C}\right)  $, the tangent
space of $\mathcal{M}^{SL}\left(  X\right)  $ at $\left(  L,D_{E}\right)  $
(\cite{Mc}).

\bigskip

It is not difficult to verify that $\mathcal{B}^{ass}\left(  M\right)
=\mathcal{B}^{SL}\left(  X\right)  \times S^{1}$ and the decomposition for
$\mathcal{M}^{ass}\left(  M\right)  $ and $\Omega_{\mathcal{M}^{ass}\left(
M\right)  }$.

\begin{center}
{\large On the moduli of coassociative cycles}
\end{center}

Suppose $C$ is a coassociative submanifold of product type in $M=X\times
S^{1}$ then

\begin{center}%
\begin{tabular}
[c]{|l|l|}\hline
Case (1): $C=L\times S^{1}$ & $L$ special Lagrangian with phase = $\pi/2\
\begin{array}
[c]{l}%
\,\\
\,
\end{array}
$\\\hline
Case (2): $C=S\times\left\{  t\right\}  $ & $S$ complex surface in $X$.$%
\begin{array}
[c]{l}%
\,\\
\,
\end{array}
$\\\hline
\end{tabular}

\end{center}

Case (1): $C=L\times S^{1}$. We can argue as in the previous situation to
obtain decompositions for $\mathcal{M}^{coa}\left(  M\right)  $ and
$\Omega_{\mathcal{M}^{coa}\left(  M\right)  }$.

\bigskip

Case (2): $C=S\times\left\{  t\right\}  $. We denote the moduli space of
complex surfaces $S$ (resp. $\left(  S,D_{E}\right)  $ with $D_{E}$ a flat
$U\left(  1\right)  $-connection on $S$) in $X$ as $\mathcal{B}^{cx}\left(
X\right)  $ (resp. $\mathcal{M}^{cx}\left(  X\right)  $). The tangent space of
$\mathcal{M}^{cx}\left(  X\right)  $ can be identified as $H^{2,0}\left(
S\right)  \oplus H^{0,1}\left(  S\right)  $ because of $K_{X}=O_{X}$. We have
natural differential forms on $\mathcal{M}^{cx}\left(  X\right)  $:

\begin{center}%
\begin{tabular}
[c]{|l|l|}\hline
holomorphic 3 form: & $\Omega_{\mathcal{M}^{cx}\left(  X\right)  }=\ \int
_{S}\phi\wedge\alpha\wedge\beta,
\begin{array}
[c]{l}%
\,\\
\,
\end{array}
$\\\hline
symplectic form: & $\omega_{\mathcal{M}^{cx}\left(  X\right)  }=\ \int
_{S}\left(  \left\langle \phi,J\eta\right\rangle +\left\langle \alpha
,J\beta\right\rangle \right)  \omega^{2}/2.
\begin{array}
[c]{l}%
\,\\
\,
\end{array}
$\\\hline
\end{tabular}

{\small Here }$\phi,\eta\in H^{2,0}\left(  S\right)  ${\small \ and }%
$\alpha,\beta,\in H^{0,1}\left(  S\right)  .$
\end{center}

Again we can verify the decompositions for $\mathcal{M}^{coa}\left(  M\right)
$ and $\Omega_{\mathcal{M}^{coa}\left(  M\right)  }$.

\section{\label{1SecFibFourier}Fibrations and Fourier transforms}

We review briefly the Strominger-Yau-Zaslow's construction \cite{SYZ} of the
conjectural mirror manifold of a Calabi-Yau threefold: If $X$ is close to a
large complex and K\"{a}hler structure limit point, then it should admit a
special Lagrangian fibration,
\[
\pi_{X}:X\rightarrow B,
\]
calibrated by $\operatorname{Im}\Omega_{X}$. Moreover there is a special
Lagrangian section, calibrated by $\operatorname{Re}\Omega_{X}$. The mirror
manifold is conjectured to be the dual torus fibration over $B$. Moreover the
fiberwise Fourier transformation is expected to play an important role in the
mirror transformation of the complex (resp. symplectic) geometry of $X$ to the
symplectic (resp. complex) geometry of its mirror manifold.

In this section we study calibrated fibrations on $G_{2}$-manifolds and
Fourier transformations along them.

\subsection{\label{1SecCalibFib}Calibrated fibrations}

\textbf{Coassociative fibrations}

Starting with a special Lagrangian fibration with a section on $X$, the
induced fibration on the $G_{2}$-manifold $M=X\times S^{1}$,
\[
\pi:M\rightarrow B,
\]
becomes a coassociative $T^{4}$-fibration with an associative section.

Recall that a general fiber in a Lagrangian fibration on any symplectic
manifold is necessarily a torus. This is not true for a general fiber in a
coassociative fibration on a $G_{2}$-manifold. For example, given any
holomorphic fibration on a Calabi-Yau threefold,
\[
p_{X}:X\rightarrow\mathbb{CP}^{1}\cong S^{2}\text{,}%
\]
its general fiber is necessarily a K3 surface. The induced fibration on the
$G_{2}$-manifold $M=X\times S^{1}$ over $S^{2}\times S^{1}$ is again a
coassociative fibration, which is now a K3 fibration.

We expect that a general fiber in a coassociative fibration is always a torus
or a K3 surface\footnote{We learn this from S.-T. Yau.}. This is indeed the
case provided that a general fiber admits an Einstein metric: The normal
bundle of a coassociative submanifold $F$ is always the bundle of self-dual
two forms on $F$. When $F$ is a fiber of $\pi$ then the normal bundle is
necessarily trivial, i.e.
\[
\Lambda_{+}^{2}\left(  F\right)  =F\times\mathbb{R}^{3}\text{.}%
\]
This implies that
\[
3\tau\left(  F\right)  +2\chi\left(  F\right)  =0,
\]
where $\tau\left(  F\right)  $ and $\chi\left(  F\right)  $ are the signature
and the Euler characteristic of $F$ respectively.

\bigskip

From physical considerations we expect that there should be a family of
$G_{2}$-metrics together with coassociative fibrations parametrized by
$\left(  t_{0},\infty\right)  $ such that the limit as $t\rightarrow\infty$
would imply the second fundamental form of each smooth fiber would go to zero.
Combining with the fact the Ricci tensor of a $G_{2}$-metric is trivial, this
suggest that every smooth fiber should admit a Einstein metric. However
characteristic numbers for $F$ saturate the Hitchin inequality for Einstein
four manifolds, this implies that $F$ is either flat or covered by a K3 surface.

\bigskip

\textbf{Associative fibrations}

Suppose $M$ has an associative fibration,
\[
\pi:M\rightarrow B.
\]
As in the previous case, we expect a general fiber to admit a metric with zero
Ricci curvature in the adiabatic limit, thus it must be a three torus $T^{3}$.

\bigskip

\textbf{Conjectural characterizations}

Calibrated fibrations on a Calabi-Yau threefold $X$ can be either special
Lagrangian $T^{3}$-fibrations or holomorphic fibrations. In the latter case a
generic fiber can be an elliptic curve $T^{2}$, an Abelian surface $T^{4}$ or
a K3 surface. Wilson \cite{Wi} studies such fibrations in terms of the ring
structure on $H^{2}\left(  X\right)  $ and the linear form,
\[
\int_{X}p_{1}\left(  X\right)  \cup\left(  \cdot\right)  :H^{2}\left(
X,\mathbb{R}\right)  \rightarrow\mathbb{R}.
\]
Notice that such linear form is trivial on $H^{3}\left(  X,\mathbb{R}\right)
$.

\bigskip

We conjecture that similar characterizations should hold true for $G_{2}%
$-manifolds. To be precise, let $\Omega_{t}$ be a family of $G_{2}$-structures
on $M$ for $t\in\left(  0,1\right)  $ such that as $t\rightarrow0 $ the volume
of $M$ goes to zero while the diameter remains constant. For small $t$, $M$
should admit a calibrated fibration, the structure of this fibration should be
determined by the Yukawa coupling $\mathcal{Y}$ and $p_{1}\left(  M\right)  $
as follows,%

\[%
\begin{tabular}
[c]{|l|l|}\hline
On $M,$ & Characterizations:$%
\begin{array}
[c]{l}%
\,\\
\,
\end{array}
$\\\hline\hline
\multicolumn{1}{|r|}{Associative $T^{3}$-fibration} & $\mathcal{Y}\left(
\Omega_{t},\Omega_{t},H^{3}\right)  \rightarrow0
\begin{array}
[c]{l}%
\,\\
\,
\end{array}
$\\\hline
\multicolumn{1}{|r|}{Coassociative $T^{4}$-fibration} & $\mathcal{Y}\left(
\Omega_{t},\Omega_{t},H^{3}\right)  \nrightarrow0$ \& $\int p_{1}\cup
\Omega_{t}\rightarrow0
\begin{array}
[c]{l}%
\,\\
\,
\end{array}
$\\\hline
\multicolumn{1}{|r|}{Coassociative K3-fibration} & $\mathcal{Y}\left(
\Omega_{t},\Omega_{t},H^{3}\right)  \nrightarrow0$ \& $\int p_{1}\cup
\Omega_{t}\nrightarrow0.
\begin{array}
[c]{l}%
\,\\
\,
\end{array}
$\\\hline
\end{tabular}
\]
It is also possible that $\mathcal{Y}\left(  \Omega_{t},\Omega_{t}%
,H^{3}\right)  $ becomes unbounded, but this should only happen when
$\tilde{M}=X\times\mathbb{R}$.

Remark: We also expect that in above situations, $\lim_{t\rightarrow0}%
\Omega_{t}$ is of infinite distance from any point in the moduli space of
$G_{2}$-metrics on $M$. On the other hand finite distance boundary points in
this moduli space (i.e. incompleteness) might be related to contractions of
(i) an ADE configuration of associative $S^{3}$'s or (ii) an associative
family of ADE configuration of $S^{2}$'s and so on. C.L.\ Wang further
conjectures that such incompleteness of the moduli space should be equivalent
to degenerating families of $G_{2}$-structures on $M\ $with both the volume
and the diameter bounded away zero and infinity uniformly.

\subsection{Analog of the SYZ conjecture}

We are going to propose an analog of the SYZ conjecture for $G_{2}$-manifolds,
namely a transformation of $G_{2}$-geometry along a calibrated fibration.
Acharya \cite{Ac} is the first to use these $T^{4}$- or $T^{3}$-fibrations to
study the mirror symmetry for $G_{2}$-manifolds from a physics point of view.

\medskip\bigskip

\bigskip

\bigskip

\textbf{The }$G_{2}$ \textbf{mirror conjecture}

There is a \textit{reasonably} large class of $G_{2}$-manifolds such that
every such $M$ satisfies the followings:

(1) \textit{Mirror pair}: There is a family of $G_{2}$-metric $g_{t}$ together
with a coassociative $T^{4}$-fibration and an associative section on $M$. As
$t$ goes to infinity the diameter and the curvature tensor of each fiber go to
zero. The moduli space $\mathcal{M}^{coa}\left(  M\right)  $ of coassociative
cycles $\left(  C,D_{E}\right)  $ in $M$ when $C$ is a fiber admits a suitable
compactification $W_{t}$ as a $G_{2}$-manifold with $\Omega_{\mathcal{M}%
^{coa}\left(  M\right)  }$ the calibrating three form (after suitable
instanton corrections from associative cycles). The relation between $M$ and
$W$ is involutive.

(2) \textit{Yukawa couplings}: The Yukawa couplings $\mathcal{Y}_{M}$ and
$\mathcal{Y}_{W}$ for $M$ and $W$ on $H_{27}^{3}\left(  M\right)  $ and
$H_{27}^{3}\left(  W\right)  $ should be naturally identified after suitable
corrections coming from associative cycles in $M$ and $W$.

(3) \textit{Coassociative geometry}: There is an identification between moduli
spaces of coassociative pairs or points in $M$ and moduli spaces of
corresponding objects in $W$.

(4) \textit{Associative geometry: }There is an identification between moduli
spaces of associative pairs or deformed Donaldson-Thomas bundles in $M$ and
moduli spaces of corresponding objects in $W$.

The mirror symmetry conjecture for $G_{2}$-manifolds can be summarized in the
following table:

\begin{center}%
\begin{tabular}
[c]{|l|}\hline
Coassociative geometry on $M$
$<$%
=%
$>$%
Associative geometry on $W,
\begin{array}
[c]{l}%
\,\\
\,
\end{array}
$\\\hline
Associative geometry on $M$
$<$%
=%
$>$%
Coassociative geometry on $W$.$%
\begin{array}
[c]{l}%
\,\\
\,
\end{array}
$\\\hline
\end{tabular}

\end{center}

\subsection{\label{Sec Trans coa}Fourier transform on coassociative $T^{4}%
$-fibrations}

In this subsection we verify various parts of the above conjecture for
\textit{flat }$G_{2}$-manifolds, i.e. $M=B\times T$ with $B=\mathbb{R}^{3}$
and $T=\mathbb{R}^{4}/\Lambda$ a flat torus. The projection map
\[
\pi:M\rightarrow B,
\]
is a coassociative $T^{4}$-fibration on $M$. The large structure limit on $M$
can be obtained by rescaling $\Lambda$. If we denote $\mathcal{M}^{coa}\left(
M\right)  $ the moduli space of coassociative cycles $\left(  C,D_{E}\right)
$ on $M$ with $\left[  C\right]  $ representing a fiber class and $D_{E}$ is a
flat $U\left(  1\right)  $-connection on $C$, then it can be naturally
identified with the total space $W=B\times T^{\ast}$ of the dual torus
fibration
\[
\pi^{\prime}:W\rightarrow B\text{.}%
\]
Moreover, under this identification, the canonical three form and four form on
the moduli space correspond to the calibrating three form and four form for
the $G_{2}$-structure on $W$:
\begin{align*}
\Omega_{W}  &  =\Omega_{\mathcal{M}^{coa}\left(  M\right)  },\\
\Theta_{W}  &  =\Theta_{\mathcal{M}^{coa}\left(  M\right)  }.
\end{align*}

We are going to perform fiberwise Fourier transformation on $M$. As usual we
denote the coordinates on $B$ (resp. $T$ and $T^{\ast}$) as $x^{1},x^{2}%
,x^{3}$ (resp. $y^{0},\cdots,y^{3}$ and $y_{0},\cdots,y_{3}$). The dual torus
$T^{\ast}$ is treated as the moduli space of flat $U\left(  1\right)  $
connections on $T$. The universal line bundle over $T\times T^{\ast}$ is
called the Poincar\'{e} bundle $P$. We normalize it so that its restriction to
both $T\times0$ and $0\times T^{\ast}$ are trivial. It has a universal
connection whose curvature two form equals
\[
\mathbb{F}=%
{\textstyle\sum_{j=0}^{3}}
dy^{j}\wedge dy_{j}\in\Omega^{2}\left(  T\times T^{\ast}\right)  \text{.}%
\]

The Fourier transformation, or the Fourier-Mukai transformation, from $M$ to
$W$ can be described roughly as follows:
\[
\mathbf{F}\left(  \bullet\right)  =p_{\ast}^{\prime}\left(  p^{\ast}\left(
\bullet\right)  \otimes P\right)  ,
\]
where $p$ and $p^{\prime}$ is defined by the following commutative diagram,%

\[%
\begin{array}
[c]{ccc}%
M\times_{B}W & \overset{p^{\prime}}{\rightarrow} & W\\
p\downarrow\, &  & \,\downarrow\pi^{\prime}\\
M & \overset{\pi}{\rightarrow} & B.
\end{array}
\]
and $p^{\ast}$ is the pull back and $p_{\ast}^{\prime}$ is the push forward
operator or equivalently integration along fibers of $p^{\prime}$ $:$
$M\times_{B}W\longrightarrow W$. We can discuss this transformation on
different levels, including differential forms, cohomology, K-theory and so
on. For differential forms, we have%
\begin{align*}
\mathbf{F}  &  :\Omega^{\ast}\left(  M\right)  \longrightarrow\Omega^{\ast
}\left(  W\right) \\
\mathbf{F}\left(  \varphi\right)   &  =\int_{M/B}p^{\ast}\left(
\varphi\right)  \wedge e^{\frac{\sqrt{-1}}{2\pi}\mathbb{F}}%
\end{align*}
where $\mathbb{F}$ is the universal curvature $2$-from on the Poincar\'{e}
bundle $P.$ Moreover $\varphi$ could be a differential form supported on a
submanifold in $M$, and the same for $\mathbf{F}\left(  \varphi\right)  $,
too. For example the above Fourier transformation preserves the calibrating
three form and four form, more precisely we have,

\begin{center}%
\[%
\begin{tabular}
[c]{|l|}\hline
$\mathbf{F}\left(  e^{\Theta_{M}}\right)  =e^{\Theta_{W}},%
\begin{array}
[c]{l}%
\,\\
\,
\end{array}
$\\\hline
$\mathbf{F}\left(  \ast e^{\Theta_{M}}\right)  =\ast e^{\Theta_{W}}.%
\begin{array}
[c]{l}%
\,\\
\,
\end{array}
$\\\hline
\end{tabular}
\]

\bigskip
\end{center}

These follow from the following integrals with orientation $y^{0123}$,
\[
\int_{M/B}e^{\Theta_{M}}\wedge e^{\mathbb{F}}=e^{\Theta_{W}}\text{ and }%
\int_{M/B}\left(  \ast e^{\Theta_{M}}\right)  \wedge e^{\mathbb{F}}=\left(
\ast e^{\Theta_{W}}\right)  .
\]

Similarly, on the level of cohomology classes, we have,%
\begin{align*}
\mathcal{F}  &  :H^{\ast}\left(  M\right)  \longrightarrow H^{\ast}\left(
W\right) \\
\mathcal{F}\left(  \varphi\right)   &  =\int_{M/B}p^{\ast}\left(  \left[
\varphi\right]  \right)  \cup ch\left(  P\right)  .
\end{align*}

We can also transform a connection $D_{E}$ on a bundle $E$ over a submanifold
$C$ in $M$ to one in $W$. To do this, we consider tensor product of the
pullback connection on the pullback bundle $p^{\ast}E\longrightarrow
C\times_{B}W$ with the universal connection on the Poincar\'{e} bundle
restricted to $C\times_{B}W$. Under good circumstances, the pushforward sheaf
under $p^{\prime}:C\times_{B}W\rightarrow W$ gives a vector bundle over a
submanifold in $W$. For instance, a flat $U\left(  1\right)  $ connection over
$T^{4}\times\left\{  p\right\}  $ in $M=$ $T^{4}\times T^{3}$ is transformed
to a point in $W=T^{4\ast}\times T^{3}$, as discussed earlier.

\bigskip

\begin{center}
{\large Transforming coassociative cycles}
\end{center}

\textbf{Transforming ASD connections on }$T^{4}$\textbf{\ fibers}

Suppose that we relax the flatness assumption on $D_{E}$ and we assume that it
is an ASD connection of a torus fiber $C=T$. The moduli space of such
coassociative cycles $\left(  C,D_{E}\right)  $ has the form,
\[
\mathcal{M}^{coa}\left(  M\right)  =B\times\mathcal{M}^{ASD}\left(  T\right)
\text{,}%
\]
where $\mathcal{M}^{ASD}\left(  T\right)  $ denotes the moduli space of ASD
connections on the flat torus $T$.

When $D_{E}$ has no flat factor, the Fourier transformation of an ASD
connection over a flat torus $T$ is studied by Schenk \cite{Sc}, Braam and van
Baal \cite{BB} and it produces another ASD connection over the dual torus
$T^{\ast}$. Moreover it is an isometry from $\mathcal{M}^{ASD}\left(
T\right)  $ to $\mathcal{M}^{ASD}\left(  T^{\ast}\right)  $. Using their
results we have,
\[
\mathbf{F}:\mathcal{M}^{coa}\left(  M\right)  \overset{\approx}{\rightarrow
}\mathcal{M}^{coa}\left(  W\right)  .
\]

Moreover it is not difficult to identify the natural three forms
$\Omega_{\mathcal{M}^{coa}}$ on these moduli spaces under $\mathbf{F}$.

\bigskip

\textbf{Describing other coassociative cycles in }$M$

Suppose $C$ is any coassociative submanifold in $M$ which is not a fiber. The
intersection of $C$ with a general fiber can not have dimension 3 (see section
\ref{Sec Coa cycle}). The intersection can not have dimension 1 neither
because otherwise the normal bundle of $C$ in $M$, which is associative, lies
inside the relative tangent bundle of the coassociative fibration, a
contradiction. Therefore the intersection of $C$ with a general fiber must
have dimension 2.

\bigskip

We assume that the coassociative submanifold $C$ in $M$ is \textit{semi-flat}
in the sense that it intersects each fiber along a flat subtorus. As subtori
in $T^{4}$, they are all translations of each other. Using the fact that the
vector cross product on $M$ preserves the normal directions of $C$, one can
show that the image of $C$ in $B\cong\mathbb{R}^{3}$ must be an affine plane.
Therefore, up to coordinate changes, we can assume that $C$ is of the form
\[
C:\left\{
\begin{array}
[c]{l}%
x^{3}=0,\\
y^{0}=B^{0}\left(  x^{1},x^{2}\right)  ,\\
y^{3}=B^{3}\left(  x^{1},x^{2}\right)  ,
\end{array}
\right.
\]
for some smooth functions $B^{0}$ and $B^{3}$ in variables $x^{1}$ and $x^{2}
$ satisfying
\[
-\frac{\partial}{\partial x^{2}}B^{3}+\frac{\partial}{\partial x^{1}}%
B^{0}=\frac{\partial}{\partial x^{1}}B^{3}+\frac{\partial}{\partial x^{2}%
}B^{0}=0.
\]

Suppose $D_{E}$ is an ASD connection on $C$ which is \textit{semi-flat }in the
sense that $D_{E}$ is invariant under translations along fiber directions, we
write
\[
D_{E}=d+a_{1}\left(  x^{1},x^{2}\right)  dx^{1}+a_{2}\left(  x^{1}%
,x^{2}\right)  dx^{2}+D_{1}\left(  x^{1},x^{2}\right)  dy^{1}+D_{2}\left(
x^{1},x^{2}\right)  dy^{2},
\]
and the ASD condition is equivalent to
\begin{align*}
\frac{\partial}{\partial x^{2}}a_{1}-\frac{\partial}{\partial x^{1}}a_{2}  &
=0,\\
-\frac{\partial}{\partial x^{2}}D_{1}+\frac{\partial}{\partial x^{1}}D_{2}  &
=\frac{\partial}{\partial x^{1}}D_{1}+\frac{\partial}{\partial x^{2}}D_{2}=0.
\end{align*}

On each $T^{4}$ fiber of $\pi:$ $M\rightarrow B$ over a point $\left(
x^{1},x^{2},0\right)  \in B$, we obtain a flat $2$-torus defined by%
\[
\left\{  \left(  y^{0},y^{1},y^{2},y^{3}\right)  \in T^{4}\mid y^{0}%
=B^{0}\left(  x^{1},x^{2}\right)  ,\;y^{3}=B^{3}\left(  x^{1},x^{2}\right)
\right\}
\]
together with a flat $U\left(  1\right)  $ connection $D_{E}=d+D_{1}\left(
x^{1},x^{2}\right)  dy^{1}+D_{2}\left(  x^{1},x^{2}\right)  dy^{2}$ on it,
because $D_{1}\left(  x^{1},x^{2}\right)  $ and $D_{2}\left(  x^{1}%
,x^{2}\right)  $ are constant on each fiber.

\bigskip

\textbf{Transforming coassociative families of subtori}

BApplying the Fourier transformation from this $T^{4}$ fiber to $T^{4\ast}$,
it is not difficult to check that the above flat connection ( resp. flat
$2$-torus ) in $T^{4}$ will be transformed to a flat $2$-torus ( resp. flat
connection ) in $T^{4\ast}$ ( see \cite{LYZ} for more details ). Putting these
together over various fibers of $\pi:M\rightarrow B$, we obtain the Fourier
transformation of $\left(  C,D_{E}\right)  $ in $M$, which is the pair
$\left(  C^{\prime},D_{E}^{\prime}\right)  $ on $W$ with%

\[
C^{\prime}:\left\{
\begin{array}
[c]{l}%
x^{3}=0,\\
y_{1}=D_{1}\left(  x^{1},x^{2}\right)  ,\\
y_{2}=D_{2}\left(  x^{1},x^{2}\right)  .
\end{array}
\right.
\]
and
\[
D_{E}^{\prime}=d+a_{1}\left(  x^{1},x^{2}\right)  dx^{1}+a_{2}\left(
x^{1},x^{2}\right)  dx^{2}+B^{0}\left(  x^{1},x^{2}\right)  dy_{0}%
+B^{3}\left(  x^{1},x^{2}\right)  dy_{3}.
\]
$\left(  C^{\prime},D_{E}^{\prime}\right)  $ is again a semi-flat
coassociative cycle in $W$ and notice that the ASD connection for $D_{E}$
(resp. coassociative condition for $C$) implies the coassociative condition
for $C$ (resp. ASD condition for $D_{E}^{\prime}$). Furthermore the Fourier
transformation induces a bijection on the semi-flat part of the moduli
spaces,
\[
\mathbf{F:}\mathcal{M}^{coa}\left(  M\right)  _{semiflat}\rightarrow
\mathcal{M}^{coa}\left(  W\right)  _{semiflat}.
\]

It is easy to check that when $C$ is horizontal, i.e. $B^{0}=B^{3}=0$, then
$\mathbf{F}$ preserves the natural three forms $\Omega_{\mathcal{M}^{coa}}$ on
these moduli spaces and we expect it continues to hold in general.

\bigskip

\begin{center}
{\large Transforming associative cycles and DT bundles}
\end{center}

\textbf{Describing associative cycles in }$M$

Suppose $A$ is any associative submanifold in $M$ then it can only intersect a
general fiber in dimension 2 or 0. In the former case, we can argue as before,
the image of $A$ in $B\cong\mathbb{R}^{3}$ must be an affine line provided
that $A$ is \textit{semi-flat}. We can write, up to coordinate changes,%

\[
A:\left\{
\begin{array}
[c]{l}%
x^{2}=x^{3}=0,\\
y^{2}=B^{2}\left(  x^{1}\right)  ,\\
y^{3}=B^{3}\left(  x^{1}\right)  .
\end{array}
\right.
\]

The associative condition on $A$ implies that
\[
\frac{\partial}{\partial x^{1}}B^{2}=\frac{\partial}{\partial x^{1}}%
B^{3}=0\text{.}%
\]
That is, both $B^{2}\left(  x^{1}\right)  $ and $B^{3}\left(  x^{1}\right)  $
are constant functions.

Suppose $D_{E}$ is a unitary connection on $A$ which is invariant under
translations along fiber directions, we can write,%

\[
D_{E}=a\left(  x^{1}\right)  dx^{1}+D^{0}\left(  x^{1}\right)  dy_{0}%
+D^{1}\left(  x^{1}\right)  dy_{1}.
\]
If $\left(  A,D_{E}\right)  $ is an associative cycle, then the flatness of
$D_{E}$ is equivalent to
\[
\frac{\partial}{\partial x^{1}}a=\frac{\partial}{\partial x^{1}}D^{0}%
=\frac{\partial}{\partial x^{1}}D^{1}=0\text{.}%
\]
That is $a,D^{0}$ and $D^{1}$ are all constant functions.

\bigskip

\textbf{Transformating associative families of tori}

As in the earlier situation, the Fourier transformation of a semi-flat cycle
$\left(  A,D_{E}\right)  $ on $M$, which is not necessarily associative, is a
pair $\left(  A^{\prime},D_{E}^{\prime}\right)  $ on $W$ with%

\[
A^{\prime}:\left\{
\begin{array}
[c]{l}%
x^{2}=0,\\
x^{3}=0,\\
y_{0}=D^{0}\left(  x^{1}\right)  ,\\
y_{1}=D^{1}\left(  x^{1}\right)  .
\end{array}
\right.
\]
and%

\[
D_{E}^{\prime}=a\left(  x^{1}\right)  dx^{1}+B^{2}\left(  x^{1}\right)
dy_{2}+B^{3}\left(  x^{1}\right)  dy_{3}%
\]

It is obvious that $\left(  C,D_{E}\right)  $ is associative if and only if
$\left(  C^{\prime},D_{E}^{\prime}\right)  $ is associative.

Remark: Recall that associative cycles in $M$ are critical points of the
Chern-Simons functional. Using arguments in \cite{LYZ}, we can show that the
Fourier transformation preserves the Chern-Simons functional of semi-flat
three cycles on $M$ and $W$ which are \textit{not} necessarily associative.

\bigskip

\textbf{Transforming associative sections to deformed DT bundles}

Any section of the (trivial) coassociative fibration $\pi$ on $M$ can be
expressed as a function $f$,
\[
f:B\rightarrow T=\frac{\mathbb{H}}{\Lambda}\text{.}%
\]
Harvey and Lawson \cite{HL} show that the associative condition is equivalent
to the following equation,
\[
-\frac{\partial f}{\partial x^{1}}i-\frac{\partial f}{\partial x^{2}}%
j-\frac{\partial f}{\partial x^{3}}k=\frac{\partial f}{\partial x^{1}}%
\times\frac{\partial f}{\partial x^{2}}\times\frac{\partial f}{\partial x^{3}%
}\text{.}%
\]
If we write $f=\,^{t}\left(  f^{0},f^{1},f^{2},f^{3}\right)  \in\mathbb{R}%
^{4}/\Lambda$, $I=\,^{t}\!\left(  1,i,j,k\right)  $ and $f_{x^{j}}^{i}%
=\frac{\partial f^{i}}{\partial x^{j}}$, then the above equation can be
rewritten as follows,
\[%
\begin{array}
[c]{cc}%
\det\left(  I,f_{x^{1}},f_{x^{2}},f_{x^{3}}\right)  _{4\times4}= & \left(
f_{x^{1}}^{1}+f_{x^{2}}^{2}+f_{x^{3}}^{3}\right)  1+(-f_{x^{1}}^{0}+f_{x^{2}%
}^{3}-f_{x^{3}}^{2})i+\\
& (-f_{x^{1}}^{3}-f_{x^{2}}^{0}+f_{x^{3}}^{1})j+(f_{x^{1}}^{2}-f_{x^{2}}%
^{1}-f_{x^{3}}^{0})k.
\end{array}
\]

If $D_{E}$ is a unitary flat connection on $A$, we can write%

\[
D_{E}=d+a_{1}dx^{1}+a_{2}dx^{2}+a_{3}dx^{3}%
\]
where $a_{k}$'s are matrix valued functions of $\left(  x^{1},x^{2}%
,x^{3}\right)  $.

As before, or from \cite{LYZ}, the Fourier transformation of $\left(
A,D_{E}\right)  $ is a unitary connection $D_{E}^{\prime}$ on the whole
manifold $W$ given by
\[
D_{E}^{\prime}=d+a_{1}dx^{1}+a_{2}dx^{2}+a_{3}dx^{3}+f^{0}dy_{0}+f^{1}%
dy_{1}+f^{2}dy_{2}+f^{3}dy_{3}.
\]
We are going to show that $D_{E}^{\prime}$ satisfies the deformed DT equation,
i.e.
\[
F_{E^{\prime}}\wedge\Theta_{W}+F_{E^{\prime}}^{3}/6=0\text{.}%
\]
First we notice that the curvature of $D_{E}^{\prime}$ is given by%

\[
\ F_{E^{\prime}}=\sum_{j=1}^{3}\sum_{k=0}^{3}f_{x^{i}}^{k}dx^{i}\wedge
dy_{k}.
\]
and $a_{j}$'s do not occur here because of the flatness of $D_{E}$. By direct
but lengthy computations, we can show that the deformed DT equation for
$D_{E}^{\prime}$ is reduced to the flatness of $D_{E}$ together with the
associativity of $A$. Therefore the Fourier transformation in this case
gives,
\[
\mathbf{F}:\mathcal{M}^{ass}\left(  M\right)  \rightarrow\mathcal{M}%
^{bdl}\left(  W\right)  \text{.}%
\]

In fact this transformation $\mathbf{F}$ preserves the Chern-Simons
functionals (see section \ref{1SecDTbdl}, \ref{1SecAssCycle}) on $M$ and $W$.
Namely suppose $\left(  A,D_{E}\right)  $ is any pair on $M$ which is not
necessarily an associative cycle, then we have
\[
\int_{A\times\left[  0,1\right]  }Tr\left[  e^{\Theta_{M}+\tilde{F}}\right]
=\left(  const\right)  \int_{W\times I}Tr\left[  e^{\Theta_{W}+\tilde
{F}^{\prime}}\right]
\]

We expect that $\mathbf{F}$ preserves the natural three forms on
$\mathcal{M}^{ass}\left(  M\right)  $ and $\mathcal{M}^{bdl}\left(  W\right)
$.

\subsection{Fourier transform on associative $T^{3}$-fibrations}

Suppose $M$ has an associative $T^{3}$-fibration,
\[
\pi:M\rightarrow B.
\]
The fiberwise Fourier transformation $\mathbf{F}$ in this case should
interchange the two types of geometry on $G_{2}$-manifolds,
\[
\text{Associative Geometry on }M\overset{\mathbf{F}}{\longleftrightarrow
}\text{Coassociative Geometry on }W\text{.}%
\]

Here $W$ is the dual three torus fibration to $\pi$ on $M$, it can be
identified as the moduli space $\mathcal{M}^{ass}\left(  M\right)  $ of
associative cycles $\left(  A,D_{E}\right)  $ on $M\;$with $D_{E}$ being a
flat $U\left(  1\right)  $-connection over a fiber $A$.

In the \textit{flat} case, we have $M=T\times\mathbb{R}^{4}$ and
$\mathcal{M}^{ass}\left(  M\right)  \cong W=T^{\ast}\times B$. Under these
identifications we have,
\begin{align*}
\Omega_{W}  &  =\Omega_{\mathcal{M}^{ass}\left(  M\right)  },\\
\Theta_{W}  &  =\Theta_{\mathcal{M}^{ass}\left(  M\right)  }.
\end{align*}

\bigskip

We can analyze semi-flat associative and coassociative cycles in $M$ and their
Fourier transformations as before. For example, semi-flat coassociative
submanifolds in $M$ are either sections or a family of affine two tori over an
affine plane in $B$. Semi-flat associative submanifolds in $M$ are either
fibers or a family of affine two tori over an affine line in $B $.

Suppose $\left(  C,D_{E}\right)  $ is a coassociative cycle in $M$ with $C$
being a section of the (trivial) associative fibration on $M$. The
coassociative condition of $C$ can be expressed in term of a differential
equation for a function $g=\left(  g^{1},g^{2},g^{3}\right)  :\mathbb{R}%
^{4}\rightarrow T=\mathbb{R}^{3}/\Lambda$ (see \cite{HL}),
\[
\nabla g^{1}i+\nabla g^{2}j+\nabla g^{3}k=\nabla g^{1}\times\nabla g^{2}%
\times\nabla g^{3}\text{.}%
\]
If we write the ASD connection as
\[
\ D_{E}=d+a_{0}dy^{0}+a_{1}dy^{1}+a_{2}dy^{2}+a_{3}dy^{3},
\]
then the Fourier transformation of $\left(  C,D_{E}\right)  $ will be the
following connection $D_{E}^{\prime}$ on $W$,
\[
D_{E}^{\prime}=d+a_{0}dy^{0}+a_{1}dy^{1}+a_{2}dy^{2}+a_{3}dy^{3}+g^{1}%
dx_{1}+g^{2}dx_{2}+g^{3}dx_{3}%
\]

Again direct but lengthy computations show that the deformed DT equation on
$D_{E}^{\prime}$ is reduced to the above coassociativity condition of $C$ and
the ASD equation for $D_{E}$.

\newpage

\section{\label{Sec Spin(7)}$Spin\left(  7\right)  $-geometry}

There are only two kinds of exceptional holonomy groups for Riemannian
manifolds, namely $G_{2}$ for dimension $7$ and $Spin\left(  7\right)  $ in
dimension $8$. Both of them play important roles in M-theory. Almost all of
our earlier discussions for $G_{2}$-geometry has a $Spin\left(  7\right)  $
analog. As a result, we will only indicate what geometric structures we have
in $Spin\left(  7\right)  $-geometry and leave out those derivations which are
the same as in their $G_{2}$-counterparts.

An eight dimensional Riemannian manifold $Z$ with $Spin\left(  7\right)
$-holonomy has a closed self-dual four form $\Theta_{Z}$ which can be
expressed as follows,
\begin{align*}
\Theta_{Z}  &  =-dy^{0123}-dx^{0123}-\left(  dx^{10}+dx^{23}\right)  \left(
dy^{10}+dy^{23}\right) \\
&  -\left(  dx^{20}+dx^{31}\right)  \left(  dy^{20}+dy^{31}\right)  -\left(
dx^{30}+dx^{12}\right)  \left(  dy^{30}+dy^{12}\right)  ,
\end{align*}
when $Z=\mathbb{R}^{8}$. As before, we can decompose differential forms on $Z$
into $Spin\left(  7\right)  $ irreducible components (see e.g. \cite{Jo}), for
example $\Lambda^{4}=\Lambda_{1}^{4}+\Lambda_{7}^{4}+\Lambda_{27}^{4}%
+\Lambda_{35}^{4}.$

\bigskip

\textbf{Holonomy reduction}

For example if $M$ is a $G_{2}$-manifold, then $Z=M\times S^{1}$ has a
$Spin\left(  7\right)  $-metric with
\[
\Theta_{Z}=\Omega_{M}\wedge dt-\Theta_{M},
\]
This corresponds to the natural inclusion $G_{2}\subset Spin\left(  7\right)
$. We can also consider the reduction of holonomy to the subgroup $Spin\left(
6\right)  \subset Spin\left(  7\right)  .$ Notice that $Spin\left(  6\right)
=SU\left(  4\right)  $, namely such a $Z$ would be a Calabi-Yau fourfold and%

\[
\Theta_{Z}=-\frac{\omega_{Z}^{2}}{2}+\operatorname{Re}\Omega_{Z}\text{,}%
\]
where $\omega_{Z}$ and $\Omega_{Z}$ are the Ricci flat K\"{a}hler form and the
holomorphic volume form on $Z$. We can compare their decompositions of
differential forms as follows,
\begin{align*}
\Lambda_{35}^{4}  &  =\Lambda_{prim}^{3,1}+\Lambda_{prim}^{1,1}+\Lambda
_{prim}^{1,3},\\
\Lambda_{27}^{4}  &  =\Lambda_{prim}^{2,2}+\Lambda^{0,2}+\Lambda^{0,4}.
\end{align*}

We can reduce the holonomy group to smaller subgroups of $Spin\left(
7\right)  $ and give the following tables,

\begin{center}%
\begin{tabular}
[c]{|l|l|}\hline
$\text{Holonomy}$ & $\text{Manifold}$\\\hline
$Spin\left(  7\right)
\begin{array}
[c]{l}%
\,\\
\,
\end{array}
$ & $Spin\left(  7\right)  \text{-manifold}$\\\hline
$Spin\left(  6\right)  =SU\left(  4\right)
\begin{array}
[c]{l}%
\,\\
\,
\end{array}
$ & $\text{Calabi-Yau 4-fold}$\\\hline
$Spin\left(  5\right)  =Sp\left(  2\right)
\begin{array}
[c]{l}%
\,\\
\,
\end{array}
$ & $\text{Hyperk\"{a}hler 4-fold}$\\\hline
$Spin\left(  4\right)  =SU\left(  2\right)  ^{2}
\begin{array}
[c]{l}%
\,\\
\,
\end{array}
$ & $K3\times K3$\\\hline
$Spin\left(  3\right)  =SU\left(  2\right)
\begin{array}
[c]{l}%
\,\\
\,
\end{array}
$ & $K3\times T^{4}$\\\hline
\end{tabular}

\end{center}

\subsection{$Spin\left(  7\right)  $-Analog of Yukawa coupling}

\textbf{An analog of the Yukawa coupling}

We consider the $\left(  3,1\right)  $-tensor $\chi_{Z}\in\Omega^{3}\left(
Z,T_{Z}\right)  $ constructed by raising an index of $\Theta_{Z}$ using the
metric tensor.

We define a symmetric quartic tensor $\mathcal{Q}$ on $H^{4}\left(
Z,\mathbb{R}\right)  $ as follows,
\begin{align*}
\mathcal{Q}  &  :%
{\textstyle\bigotimes^{4}}
\Omega^{4}\left(  Z,\mathbb{R}\right)  \rightarrow\mathbb{R}\\
\mathcal{Q}\left(  \phi_{1},\phi_{2},\phi_{3},\phi_{4}\right)   &  =\int
_{Z}\Theta_{Z}\left(  \hat{\phi}_{1},\hat{\phi}_{2},\hat{\phi}_{3},\hat{\phi
}_{4}\right)  \wedge\Theta_{Z}.
\end{align*}
Here $\hat{\phi}=\ast\left(  \phi\wedge\chi\right)  \in\Omega^{1}\left(
Z,T_{Z}\right)  $. Note that $\hat{\phi}$ is zero if $\phi\in\Omega_{27}%
^{4}\left(  Z,\mathbb{R}\right)  $. We define the \textit{Yukawa coupling}
$\mathcal{Y}$ to be the restriction of $\mathcal{Q}$ to $H_{35}^{4}\left(
Z,\mathbb{R}\right)  $. We can also define a cubic form, a quadratic form and
a linear form on $H_{35}^{4}\left(  Z,\mathbb{R}\right)  $ by evaluating
$\mathcal{Q}$ on one, two and three $\Theta_{Z}$ before restricting it to
$H_{35}^{4}\left(  Z,\mathbb{R}\right)  $.

\subsection{Deformed DT bundles and Cayley cycles}

\textbf{Deformed DT bundles}

As in the $G_{2}$ case, we define a \textit{deformed DT connection} to be a
connection $D_{E}$ over $Z$ whose curvature tensor $F_{E}$ satisfies
\[
\ast F_{E}+\Theta_{Z}\wedge F_{E}+F_{E}^{3}/6=0,
\]
or equivalently,
\[
\left[  e^{\Theta_{Z}+F_{E}+\ast F_{E}}\right]  ^{\left[  6\right]  }=0.
\]

Its moduli space $\mathcal{M}^{bdl}\left(  Z\right)  $ carries a natural four
form $\Theta_{\mathcal{M}^{bdl}\left(  Z\right)  }$ defined by,
\[
\Theta_{\mathcal{M}^{bdl}\left(  Z\right)  }\left(  \alpha,\beta,\gamma
,\delta\right)  =\int_{Z}Tr_{E}\left[  \alpha\wedge\beta\wedge\gamma
\wedge\delta\right]  _{skew}\wedge e^{\Theta_{Z}+F_{E}}\text{.}%
\]
with $\alpha,\beta,\gamma,\delta\in H^{1}\left(  Z,ad\left(  E\right)
\right)  , $ the first cohomology group of the elliptic complex,
\[
0\rightarrow\Omega^{0}\left(  Z,ad\left(  E\right)  \right)  \overset{D_{E}%
}{\rightarrow}\Omega^{1}\left(  Z,ad\left(  E\right)  \right)  \overset
{\pi_{7}\circ D_{E}}{\rightarrow}\Omega_{7}^{2}\left(  Z,ad\left(  E\right)
\right)  \rightarrow0.
\]

\bigskip

\textbf{Cayley cycles}

A submanifold $C$ in $Z$ of dimension four is called a \textit{Cayley}
submanifold if it is calibrated by $\Theta_{Z}$, i.e. The restriction on
$\Theta_{Z}$ to $C$ equals the volume form on $C$ with respect to the induced
metric (\cite{HL}). A \textit{Cayley cycle} is defined to be any pair $\left(
C,D_{E}\right)  $ with $C$ a Cayley submanifold in $Z$ and $D_{E} $ is an ASD
connection on a bundle $E$ over $C$. The moduli space $\mathcal{M}%
^{Cay}\left(  Z\right)  $ of Cayley cycles in $Z$ also have a natural four
form $\Theta_{\mathcal{M}^{Cay}\left(  Z\right)  }$ defined by
\[
\Theta_{\mathcal{M}^{Cay}\left(  Z\right)  }=\left\{
\begin{array}
[c]{c}%
\int_{C}Tr\alpha\wedge\beta\wedge\gamma\wedge\delta\medskip\,\\
-\int_{C}\left\langle \phi,Tr(\widehat{\alpha\wedge\beta)}\cdot\eta
\right\rangle \Theta_{Z}\medskip\,\\
-\int_{C}\det\left(  \phi,\eta,\xi,\zeta\right)  \Theta_{Z}%
\end{array}
\right.
\]
where $\alpha,\beta,\gamma,\delta\in\Omega^{1}\left(  C,ad\left(  E\right)
\right)  $ and $\phi,\eta,\xi,\zeta\in Ker\mathbf{D}$ are harmonic spinors and%

\[
\widehat{\alpha\wedge\beta}=\alpha\wedge\beta+\ast(\alpha\wedge\beta
)\;\in\Omega_{+}^{2}\left(  C,ad\left(  E\right)  \right)  =\operatorname{Im}%
(\mathbb{H})\otimes ad\left(  E\right)
\]

We leave it to our readers to verify that the Fourier transformation along
Cayley $T^{4}$-fibration on a flat $Spin\left(  7\right)  $-manifold $Z$ will
transform a Cayley section cycle $\left(  C,D_{E}\right)  $ on $Z$ to a
deformed DT bundle over the dual torus fibration. Transformations of certain
non-section Cayley cycles on $Z$ can also be analyzed as in the $G_{2}$ case.

\bigskip

Acknowledgments: The second author is partially supported by NSF/DMS-0103355
and he thanks R. Bryant, R. Thomas, C.L. Wang, S.-T. Yau and E. Zaslow for
insightful and valuable discussions. He is grateful to S.-T. Yau who brought
us to the attention of this subject.

\bigskip

{\scriptsize Addresses: }

{\scriptsize Naichung Conan Leung (leung@ims.cuhk.edu.hk)}

{\scriptsize Institute of Mathematical Sciences and Department of Mathematics,
The Chinese University of Hong Kong, Hong Kong}

{\scriptsize Jae-Hyouk Lee (jhlee@math.wustl.edu)}

{\scriptsize Department of Mathematics, Washington University in St. Louis,
U.S.A.}

\end{document}